\documentclass[12pt]{amsart}
\usepackage{epsf,amscd,amssymb}
 
\newtheorem{thm}{Theorem}[section]
\newtheorem{lem}[thm]{Lemma}
\newtheorem{cor}[thm]{Corollary}
\newtheorem{prop}[thm]{Proposition}

\newtheorem{thmA}{Theorem A}

\newtheorem{thmB}{Theorem B}

\newtheorem{thmC}{Theorem C}

\theoremstyle{remark} 
\newtheorem{rem}[thm]{Remark}
\newtheorem{hyp}[thm]{Hypothesis}

\theoremstyle{definition}
\newtheorem{defn}[thm]{Definition}

\newcommand\bR{{\mathbb{R}}}
\newcommand\bC{{\mathbb C}}
\newcommand\bE{{\mathbb E}}
\newcommand\bZ{{\mathbb Z}}

\newcommand\bH{{\mathbb{H}}}

\newcommand\PO{{\rm PO}}
\newcommand\PSL{{\rm PSL}}

\newcommand\dev{{\bf dev}}
\newcommand\SI{{\bf S}}

\newcommand\ra{\rightarrow}

\newcommand\emp{\emptyset}
\newcommand\eps{\epsilon}

\newcommand\vth{\vartheta}

\newcommand\ovl{\overline}

\newcommand\calR{{\mathcal R}}
\newcommand\calS{{\mathcal S}}
\newcommand\calT{{\mathcal T}}
\newcommand\calC{{\mathcal C}}

\newcommand\hyps{\mathbb{H}^3} 
\newcommand\hypp{\mathbb{H}^2}
\newcommand\tri{\triangle}

\newcommand\dist{\mathbb{d}}
\newcommand\Isom{\rm{Isom}}
\newcommand\Area{\rm{Area}}
\newcommand\convh{\rm{convh}}

\newcommand\sing{\rm sing}
\newcommand\core{\mathcal C}
\newcommand\cl{\bf c}

\newcommand\cat{{\rm{CAT}}}

\newcommand\rmc{\hbox{\rm :}}

\begin{document}

\title[Drilling hyperbolic $3$-manifolds]
{Drilling cores of hyperbolic $3$-manifolds to prove tameness }
\author{Suhyoung Choi}
\address{Department of Mathematics \\
Korea Advanced Institute of Science and Technology\\  
305--701 Daejeon, South Korea}
\email{schoi@math.kaist.ac.kr}
\date{October 15, 2004 (Preliminary version 2)} 
\thanks{The author gratefully acknowledges  
support from Korea Research Foundation Grant (KRF-2002-070-C00010).}

\begin{abstract}
We sketch a proof of the fact that a hyperbolic $3$-manifold $M$ 
with finitely generated fundamental group and with no parabolics 
are topologically tame. This proves the Marden's conjecture. 
Our approach is to form 
an exhaustion $M_i$ of $M$ and modify the 
boundary to make them $2$-convex. 
We use the induced path-metric, which  
makes the submanifold $M_i$ $\delta$-hyperbolic
and with Margulis constants independent of $i$.
By taking the convex hull in the cover of $M_i$ 
corresponding the core, we show that 
there exists an exiting sequence of surfaces $\Sigma_i$.
We drill out the covers of $M_i$ by a core $\core$ again 
to make it $\delta$-hyperbolic. 
Then the boundary of the convex hull of $\Sigma_i$ is shown 
to meet the core. By the compactness argument of Souto, 
we show that infinitely many of $\Sigma_i$ are homotopic in $M - \core^o$. 
\end{abstract}  

\maketitle

\tableofcontents


Recently, Agol initiated a really interesting 
approach to proving Marden's conjecture
by drilling out closed geodesics. 
We started on the same general direction nearly the same time
as Agol. In this paper, we use truncation and drill out compact cores.
We do use the Agol's idea of using covering spaces and 
taking convex hulls of the cores, a method originating 
from Freedman. Our approach is somewhat different
in that we do not use end-reductions and pinched Riemannian 
hyperbolic metrics; however, 
we use the incomplete hyperbolic metric itself.
The hard geometric analysis and geometric convergence techniques can be avoid 
using the techniques of this paper. 
Except for developing a rather elementary theory of 
deforming boundary to make the submanifolds of codimension $0$ 
Gromov hyperbolic, we do not need any other highly developed 
techniques.
Note also that there is a recent paper by Calegari and Gabai
\cite{CalGab} using modified least area surfaces and closed geodesics.
The work here is independently developed from their line of 
ideas.

In this paper, we let $M$ be a hyperbolic $3$-manifold 
with a core homeomorphic to a compression body. 
Suppose that $M$ has a finitely generated fundamental group
and the holonomy is purely loxodromic
and has ends $E, E_1, \dots, E_n$. 
Let $F_1, \dots, F_n$ be the incompressible surfaces in 
neighborhoods of the ends $E_1, \dots, E_n$.  
Let $N(E)$ be a neighborhood of an end $E$ with no incompressible surface 
associated. 

The Marden's conjecture states that a hyperbolic $3$-manifold 
with a finitely generated fundamental group is homeomorphic to the interior of 
a compact $3$-manifold. It will be sufficient to prove for 
the above $M$ to prove Marden's conjecture. 

The cases when the group contains parabolic elements are left out. 
But we believe that the arguments in this paper essentially go through. 

\begin{thmA}
Let $M$ be as above with ends $E, E_1, \dots, E_n$,  
and $\core$ be a compact core of $M$. Then 
$E$ has an exiting sequence of surfaces of genus equal to 
that of the boundary component of $\partial \core$ corresponding to $E$.  
\end{thmA} 

The following implies the proof of Marden's conjecture
since all finitely generated hyperbolic manifolds 
without parabolic elements
has ends isometric to the manifolds $M$ as described here.

\begin{thmB}
$M$ is tame; that is, 
$M$ is homeomorphic to the interior of a compact manifold. 
\end{thmB}

This paper has three parts: 
In Part 1, 
let $M$ be a codimension $0$ submanifold of a hyperbolic 
$3$-manifold $N$ of infinite volume with certain nice 
boundary conditions. $M$ is locally finitely 
triangulated. Suppose that $M$ 
is $2$-convex in $N$ in the sense that any tetrahedron $T$ in $N$ with 
three of its side in $M$ must be inside $M$. 
Now let $L$ be another compact codimension $0$-submanifold $M$ 
so that $\partial L$ is incompressible in $M$ with a number 
of closed geodesics $\cl_1, \dots, \cl_n$ removed.
We suppose that $L$ is finitely triangulated. 
Given $\eps >0$, 
we show that $\partial L$ can be isotopied to a hyperbolically
triangulated surface so that it bounds in $M$
$2$-convex submanifold whose $\eps$-neighborhood contains 
$\cl_1, \dots, \cl_n$. 
The isotopy techniques will be PL-type arguments and deforming by 
crescents. An important point to be used in the proof is 
that the crescents avoid closed geodesics and geodesic laminations. 
Thus, the isotopy does not pass through the closed geodesics. 

Part 2 is as follows: 
A general hyperbolic manifold is a manifold with boundary modeled 
on subdomains in the hyperbolic space. 
A general hyperbolic manifold 
is $2$-convex if every isometry from a tetrahedron with 
an interior of its side removed extends to the tetrahedron itself. 
We show that a $2$-convex general hyperbolic manifold is 
Gromov hyperbolic. The proof is based on 
the analysis of the geometry of the vertices of the boundary 
required by the $2$-convexity.
We will also define h-surfaces as 
a triangulated surface where each triangle gets mapped to 
geodesic triangles and the sum of the induced angles at each 
vertex is always greater than or equal to $\pi$. We 
show the area bound of such surfaces. Finally, we 
show that the boundary of the convex hull of a core in 
a general hyperbolic manifold with finitely generated fundamental
group can be deformed to a nearby h-surface. This 
follows from local analysis of geometry.

In Part 3, we will give the proof of Theorems A and B
using the results of Part 1 and 2. The outline is 
given in the abstract and in Section 11. The proof 
itself is rather short spanning 9-10 pages only.

\part{$2$-convex hulls of submanifolds of 
hyperbolic manifolds}

\section{Introduction to Part 1}


The purpose of this paper is to deform a submanifold 
of such a manifold. Suppose $M$ is a submanifold of 
a complete hyperbolic manifold $N$ with finitely generated 
fundamental group. We suppose that $N$ has incompressible 
ends except for one end and is homotopy equivalent to 
a compression body. Suppose that $M$ contains 
a number of closed geodesics, and that $\partial M$ 
is incompressible in $N$ with these geodesics removed. 

We modify $M$ so that $M$ becomes Gromov hyperbolic and 
its $\eps$-neighborhood still contains the closed geodesics.




We will be working in a more general setting.
A {\em general hyperbolic manifold} is a Riemannian manifold $M$ with 
corner and a geodesic metric that admits
a geodesic triangulation so that each
$3$-simplex is isometric with a hyperbolic 
one. $M$ admits a local isometry, so-called developing map, 
$\dev: \tilde M \ra \hyps$ equivariant with respect to
a homomorphism $h: \pi_1(M) \ra \PSL(2, \bC)$. 
The pair $(\dev, h)$ is only determined up to action 
\[(\dev, h) \mapsto (g\circ \dev, g\circ h(\dot) \circ g^{-1})
\hbox{ by }  
g \in\PSL(2, \bC). \]

We remark that in Thurston's notes \cite{Thnote} 
a locally convex general hyperbolic manifold 
is shown to be covered by a domain in $\hyps$ or, 
equivalently, it can be 
extended to a complete hyperbolic manifold. 



A general hyperbolic manifold $M$ is {\em $2$-convex} if 
given a hyperbolic $3$-simplex $T$, 
a local-isometry $f: T - F^o \ra M^o$ for a face $F$ of $T$
extends to an isometry $T \ra M^o$. 
(See \cite{psconv} for more details. Actually 
projective version applies here by the Klein 
model of the hyperbolic $3$-space.)

By a {\em totally geodesic hypersurface}, we mean the union of 
components of the inverse image under a developing
map of a totally geodesic plane in $\hyps$. 
A {\em local totally geodesic hypersurface} is 
an open neighborhood of a point in the hypersurface. 

For a point, a {\em local half-space} is the closure 
in the ball of the component of an open ball around it with 
a totally geodesic hypersurface passing through it. 
The local totally geodesic hypersurface intersected with
the local-half space is said to be the {\em side} 
of the local half-space. A local half-space with 
its side removed is said to be an {\em open} local
half-space.


A surface $f: S \ra M$ is said to be {\em triangulated} 
if $S$ is triangulated and each triangle is mapped 
to a totally geodesic triangle in $M$. 
(We will generalize this notion a bit.)
An interior vertex of $f$ is an {\em s-vertex} if 
every open local half-space associated with 
the vertex does not contain the local image 
of $f$ with the vertex removed. 
A {\em strict s-vertex} is an s-vertex where 
every associated closed local half-space does not 
contain the local image.
An interior vertex of $f$ is a {\em convex vertex} if 
a local open half-space associated with
the vertex contains the local image of 
$f$ removed with the vertex. 
(Actually, we can see these definitions 
better by looking at the unit 
tangent bundle at the vertex:
Let $U$ be the unit tangent bundle at the vertex. 
Then the image of $f$ corresponds to a path 
in $U$.)
An interior vertex of $f$ is an {\em h-vertex} if 
the sum of angles of the triangles around a vertex 
is $\geq 2\pi$. 
An s-vertex is an h-vertex by Lemma \ref{lem:lengthhemisphere}.

$f:S \ra M$ is an s-map if each vertex is an s-vertex,
and $f$ is an h-map if each vertex is an h-one.
An s-map is an h-map but not conversely in general.
For an imbedding $f$ and an orientation, a convex vertex is 
said to be a {\em concave vertex} if the local half-space is 
in the exterior direction. Otherwise, the convex vertex is 
a {\em convex vertex}.
We also know that if the boundary of 
a general hyperbolic submanifold $N$ of $M$ is s-imbedded,
then $N$ is $2$-convex (see Proposition {prop:s-bd2-conv}).



\begin{thmC} 
Let $N$ be an orientable $2$-convex general hyperbolic 3-manifold. 
Let $M$ be a compact codimenion-zero submanifold of $N$ with boundary
$\partial M$. Suppose that each component of 
$\partial M$ is incompressible in $M$ if 
we remove a finite number of the image $\cl_1, 
\dots, \cl_n$ of 
the closed geodesics in the interior $M^o$ of $M$. 
Then for arbitrarily given small $\eps > 0$,
we can isotopy $M$ to a homeomorphic general hyperbolic 
$3$-manifold $M'$ in $N$ so that $M'$ is 2-convex
and an $\eps$-neighborhood of 
$M'$ contains the collection of closed geodesics.
\end{thmC}

We say that $M'$ obtained from $M$ by the above 
process is a {\em $2$-convex hull} of $M$. 
Although $M'$ is not necessarily a subset of $M$,
the curves $\cl_1, \dots, \cl_n$ is a subset of an $\eps$-neighborhood 
of $M'$ and hence we have certain amount of control. 
(Here the geodesics are allowed to self-intersect.)





In section one we review some hyperbolic manifold theory. 
We also introduce s-vertices and relationship with 
$2$-convexity.

In section 2, we introduce crescents and discuss their properties.
In section 3, we will study the first step of such 
an isotopy called the crescent move. 
In section 4, we discuss how to modify the infinite
pleating lines that can result from the crescent move 
to a triangulation. 
In section 5, we will combine the techniques of 
sections 2 and 3 to produce our move. 

In section 2, we introduce so-called crescents. 
We take the inverse image $\tilde \Sigma$ of the surface 
$\Sigma$. We assume that $\Sigma$ is incompressible 
in the ambient $2$-convex general hyperbolic 
manifold with a number of geodesics 
$\cl_1, \dots, \cl_n$ removed.
A crescent is a connected domain bounded by a totally geodesic hypersurface 
and an open surface in $\tilde \Sigma$. 
The portion of boundary in the totally geodesic hypersurface 
is said to be the $I$-part and the portion in 
$\tilde \Sigma$ is said to be the $\alpha$-part. 
A crescent may contain another crescents and so on. 
The folding number of 
a crescent is the maximum intersection number of 
the generic path from the outer part in the surface 
to the innermost component of the crescent with 
the surface removed. We show that for given $\Sigma$,
the folding number is bounded above.

A highest-level crescent is an innermost one that is contained 
in a crescent with highest folding number. 
We show that a highest-level crescent is always contained 
in an innermost crescent; i.e., so called the secondary 
highest-level crescent. In a secondary highest level crescent, 
the closure of the $\alpha$-part and the $I$-part are isotopic.
We also show that the secondary highest-level crescents meet nicely
extending their $\alpha$-parts in $\tilde \Sigma$, following 
\cite{psconv}.

In section 3, we introduce the crescent-isotopy theory. 
This is a theory to isotopy a surface in a general hyperbolic 
manifold so that all of its vertices become s-vertices. 

We form the union of secondary highest-level crescents
and can isotopy the union of 
their $\alpha$-parts to the complement $I$ in 
the boundary of their union. This is essentially 
the crescent move. 

However, there might be some parts of $\tilde \Sigma$ 
meeting $I$ tangentially from above. We need 
to first push these parts upward first
using truncations.

Also, after the move, there might be pleated parts 
which are not triangulated. We have to perturb these parts 
to triangulated parts. This forms section 4. 

In section 5, we gather our results to produce the move 
to obtain s-imbedded surface isotopic to $\Sigma$.
Our steps are as follows: 
We take the highest folding number and take all 
outer secondary highest-level crescents, do some truncations, and do 
the crescent isotopy. Then we perturb. 
Next, we take all inner highest-level crescents, do some truncations,  
and do crescent move and perturb. 
Now the highest folding number decreases by one, 
and we do the next step of the induction
until we have no crescents any more. 
In this case, all the vertices are s-vertices. 

We apply our result to a codimension-zero submanifold $M$
which contains a number of closed geodesics $\cl_1, \dots, \cl_n$. 
We assume that the boundary incompressible in the ambient manifold 
with the geodesics removed. This will prove Theorem A.

\section{Preliminary} 


In this section, we review the hyperbolic space 
and the Kleinian groups briefly. We 
discuss the relationship between the $2$-convexity of 
general hyperbolic manifolds 
and the s-vertices of the boundary components.

\subsection{Hyperbolic manifolds}

The hyperbolic $n$-space is a complete Riemannian metric 
space $(\bH^n, d)$ of constant curvature equal to $-1$. 
We will be concerned about hyperbolic plane and 
hyperbolic spaces, i.e., $n=2,3$, in this paper. 

The upper half space model for $\hypp$ is the pair 
\[(U^2, \PSL(2, \bR) \cup \overline{\PSL}(2, \bR))\]
where $U^2$ is the upper half space.

The Klein model of $\hypp$ is the pair $(B^2, \PO(1, 2))$
where $B^2$ is the unit disk and $\PO(1,2)$ is 
the group of projective transformations acting on $B^2$. 

A {\em Fuchsian} group is a discrete subgroup of 
the group of isometries of $\PO(1, 2)$ of the group 
$\Isom(\hypp)$ of isometries of $\hypp$.

There are many models of the hyperbolic $3$-space.
We shall use the upper-half space model or 
Klein model, whichever is more convenient at the time. 

The upper half-space model consists of 
the upper half-space $U$ of $\bR^3$ and 
the group of isometries are identified as 
the group of similarities of $\bR^3$ preserving $U$, 
which is identified as the union of $\PSL(2, \bC)$ and 
its conjugate $\overline{\PSL}(2, \bC)$. 

The Klein model consists of the unit ball 
in $\bR^3$ and the group of isometries are 
identified with the group of projective transformations 
preserving the unit ball, which is identified 
as $\PO(1, 3)$. 

A {\em Kleinian} group is a discrete subgroup of the group of 
isometries $\PO(1, 3)$ of the group $\Isom(\hyps)$ 
of isometries of $\hyps$. 

For the purposes of this paper, it is more convenient to
use the Klein model.

A {\em parabolic} element $\gamma$ of a Kleinian group is 
a nonidentity element such that $\dist(\gamma(x), x)$, $x \in \hyps$, has 
no lower bound other than $0$.
A {\em loxodromic} element of a Kleinian group is 
an isometry with a unique invariant axis. 
A {\em hyperbolic} element is a loxodromic 
one with invariant hyperplanes. 

For Fuchsian groups, a similar terminology holds. 

In this paper, we will restrict our Kleinian groups 
to be torsion-free and have no parabolic elements
and all elements are orientation-preserving.

\subsection{s-vertices}
We now classify the vertices of a triangulated map $f:S \ra M$. 
We do not yet require the general position property of $f$ 
but identify the vertex with its image. 

By a straight geodesic in a general hyperbolic manifold we 
mean a geodesic that maps to geodesics in $\hyps$
under the developing maps.

\begin{lem}\label{lem:vertex}
Let $f:S \ra M$ be a triangulated map. 
\begin{itemize}
\item An interior vertex of $S$ 
is either a convex-vertex, a concave vertex, or an s-vertex.
\item An s-vertex which is not strict one is contained in 
a unique local half-space with side containing 
the vertex and has to be one of the following: 
\begin{itemize}
\item a vertex with a totally geodesic local image.
\item a vertex on a edge of two totally geodesic planes 
where $f$ locally maps into one sides of each.
\item A vertex where a pair of triangles or 
edges has geodesic segments extending each other.
\end{itemize}
\item If $f$ is a general position map, then 
an s-vertex is a strict s-vertex.
\end{itemize}
\end{lem}
\begin{proof} 
Straightforward.
\end{proof}



\begin{lem}\label{lem:deforming} 
A s-vertex can be deformed to a strict s-vertex by 
an arbitrarily small amount by pushing if necessary the vertex  
from the boundary of the closed local half-space
containing the local image of $f$ in the direction of 
the open half-space. 
\end{lem}
\begin{proof} 
If the s-vertex is a strict one, then we leave it alone.
If the s-vertex is not a strict one, a closed local half-space
contains the local image of $f$.
Let $U$ be the unit tangent bundle at the vertex.
A closed hemisphere $H$ contains the path corresponding to 
the local image of $f$. 
We find an antipodal pair of points $v$ and $-v$ 
on $U$ corresponding to the local image of $f$.
Let $l$ be the spherical geodesic on $U$ connecting $v$ and 
$-v$ so that $l$ separates a nonantipodal pair of points $x$ and 
$y$ on the path of $f$. Let $\ovl{xy}$ be the minor 
geodesic connecting $x$ and $y$. Then $\ovl{xy}$ meets $l$. 
If we push our vertex in the open hemisphere direction, 
then $l$ becomes a geodesic segment of length $> \pi$ and
it meets $\ovl{x'y'}$ for $x'$ and $y'$ corresponding 
to $x$ and $y$ respectively. The conclusion follows from 
Lemma \ref{lem:cross}.
\end{proof}

\begin{lem}\label{lem:cross} 
Suppose that $f:\Sigma \ra M$ is a triangulated map. 
Let $v$ be a vertex of $f$ and $f':U_v \ra U^M_v$ be 
the induced map from the link of $v$ in $\Sigma$
to that of $v$ in $M$. Suppose that there exists 
a segment $l$ of length $> \pi$ in $U^M_v$ with 
endpoints in the image of $f'$ separating 
two points in the image of $f'$ so that the minor 
arc $\ovl{xy}$ meets $l$ transversely.
Then $v$ is an s-vertex.  
\end{lem}
\begin{proof} 
Let $z, t$ be the endpoints of $l$ and 
$x, y$ the separated points.
As above, let $H_1$, $H_2$, and $H_3$
be hemispheres with boundary geodesics 
containing $\ovl{xy}, \ovl{yz},$ and 
$\ovl{zx}$. Then it is clearly that any hemisphere 
containing $x,y,z$ is a subset of 
the union $H_1 \cup H_2 \cup H_3$. Therefore, 
the image of $f'$ is not contained in any hemisphere and 
$v$ is an s-vertex. 
\end{proof}

Given an oriented surface, a convex vertex is 
either a {\em convex vertex} or a {\em concave vertex} 
depending on whether the supporting local half-space 
is in the outer normal direction or in the inner normal 
direction.

\begin{prop}\label{prop:alts-vertex} 
A vertex of an oriented imbedded triangulated surface 
is either an s-vertex or a convex vertex or 
a concave vertex.
\end{prop}
\begin{proof}
Straightforward.
\end{proof}

In this paper, we will consider only metrically 
complete submanifolds, i.e., locally compact ones. 

\begin{prop}\label{prop:s-bd2-conv} 
A general hyperbolic manifold
$M$ is $2$-convex if and only if each 
vertex of $\partial M$ is a convex vertex or an s-vertex.
\end{prop}
\begin{proof}
Suppose $M$ is $2$-convex.
If a vertex $x$ of $\partial M$ is a concave, 
we can find a local half-open space in $M$ with its side 
passing through $x$. The side meets $\partial M$ 
only at $x$. From this, we can find a $3$-simplex 
inside with a face in the side. This contradicts 
$2$-convexity of $M$. 

Conversely, suppose that $\partial M$ has only convex vertices or s-vertices.
Let $f:T - F^o\ra M^o$ be a local-isometry from 
a $3$-simplex $T$ and a face $F$ of $T$. 
We may lift this map 
to $\tilde f: T - F^o \ra \tilde M^o$ where 
$\tilde M$ is the universal cover of $M$.
Then $\tilde f$ is an imbedding. 
Since $\tilde M$ is metrically complete, 
$\tilde f$ extends to $\tilde f': T \ra \tilde M$. 

Suppose that $f$ does not extend to $f':T \ra M^o$. 
This implies that $\tilde f'(F)$ meets $\partial \tilde M$. 
$\tilde f'(\partial F)$ does not meet $\partial \tilde M$. 
The subset $K = \partial \tilde M \cap \tilde f'(F)$ 
has a vertex $x$ of $\partial \tilde M$ 
which is an extreme point of the convex 
hull of $K$ in the image of $F$. 
There exists a supporting line $l$ at $x$ for $K$. 
We can tilt $\tilde f'(T)$ by $l$ a bit and 
the new $3$-simplex meets $\partial \tilde M$ at $x$ 
only. This implies that $x$ is not an s-vertex but 
a concave vertex, a contradiction.
\end{proof}




\section{2-convex hulls of hyperbolic manifolds} 


Let $M$ be a metrically complete $2$-convex 
general hyperbolic manifold
from now on and $\tilde M$ its universal cover. 
Let $\Gamma$ denote the deck transformation group of 
$\tilde M \ra M$. 

Let $\Sigma$ be a properly imbedded compact subsurface of 
an orientable general hyperbolic manifold $M$.
$\Sigma$ may have more than one components.  
We denote by $\tilde \Sigma$ the inverse image of 
$\Sigma$ in the universal cover $\tilde M$ of $M$. 
($\tilde \Sigma$ is not connected in general and 
components may not be universal covers of $\Sigma$.)
We assume that the triangulated $\tilde M$ is in general position
and so is $\tilde \Sigma$ under the developing maps. 

For each component $\Sigma_0$ of $\Sigma$ and 
a component $\Sigma_0'$ of $\tilde \Sigma$
mapping to $\Sigma_0$, there exists a subgroup 
$\Gamma_{\Sigma_0'}$ acting on $\Sigma_0'$
so that the quotient space is isometric to $\Sigma_0$.

\begin{hyp}\label{hyp:incompressible}
We will now assume that $\Sigma$ is incompressible in 
$M$ with a number of straight closed geodesics
$\cl_1, \dots, \cl_n$ removed. 
\end{hyp}

First, we introduce crescents 
for $\tilde \Sigma$ which is the inverse image 
of a surface $\Sigma$ in a $2$-convex general hyperbolic manifold.
We define the folding number of crescents and 
show that they are bounded above. 

We define the highest level 
crescents, i.e., the innermost crescents in the 
crescent with the highest folding number
incurring the highest folding number.
We show that closed geodesics avoid the interior of crescents. 
Given a highest level crescent, we show that
there is an innermost crescent that has a connected 
$I$-part to which the closure of the $\alpha$-part is isotopic
in the crescent by the incompressibility of $\Sigma$.
This is the secondary highest-level crescents.
We show that the secondary highest-level crescent is 
homeomorphic to its $I$-part times the unit interval.

Next, we show that if two highest-level crescents meet each other 
in their $I$-parts tangentially, then they both are included 
in a bigger secondary highest-level crescent.

Furthermore, if two secondary highest-level crescents meet
in their interiors, then they meet nicely extending their 
$\alpha$-parts. This is the so-called transversal intersection 
of two crescents.


\begin{defn}\label{defn:crescent} 
A {\em crescent} $\mathcal R$ for $\tilde \Sigma$ is 
\begin{itemize}
\item a connected domain in $\tilde M$ which is a closure of 
a connected open domain in $\tilde M$, 
\item so that its boundary is a disjoint union of 
a (connected) open domain in $\tilde \Sigma$ and 
the closed subset that is 
the disjoint union of totally geodesic $2$-dimensional domains 
in $\tilde M$ that develops into 
a common totally geodesic hypersurface in $\hyps$ under $\dev$. 
\end{itemize}
We denote by $\alpha_{\mathcal R}$ the domain in $\tilde \Sigma$ 
and $I_{\mathcal R}$ the union of totally geodesic domain. 
To make the definition canonical, we require 
$I_{\mathcal R}$ to be the maximal totally geodesic set in 
the boundary of $\mathcal R$.
We say that $I_\calR$ and $\alpha_\calR$ the {\em $I$-part} 
and the {\em $\alpha$-part} of $\calR$. 


As usual $\Sigma$ is oriented so that there are outer and 
inner directions to normal vectors. 

The subset $\tilde \Sigma \cap \mathcal R$ may have more than 
one components. For each component of ${\mathcal R} - \tilde \Sigma$,
we can assign a folding number which is the minimal generic 
intersection number of that a path from $\alpha_{\mathcal R}$
meeting $\tilde \Sigma$ to reach to the component. 
The maximum value of 
the folding number of the components is the folding number of 
the crescent $\mathcal R$.
\end{defn}

The {\em $I$-part hypersurface} is the maximal 
totally geodesic hypersurface in 
containing the $I$-part of the crescent.

A noncompact domain will be called a {\em crescent} if 
it is bounded by a (connected) domain in $\tilde \Sigma$ and 
the union of totally geodesic domains in 
developing into a common totally geodesic plane 
in $\hyps$ and is a geometric limit of 
compact crescents. Again the $I$-part is the maximal totally geodesic 
subset of the boundary of the crescent.
The $\alpha$-part is the complement in the boundary 
of the crescent. 

A priori, a crescent may have an infinite folding number. 
However, we will soon show that the folding number is finite. 

A crescent is an {\em outer} one if the adjacent part to 
the $\alpha$-part is in the outer normal direction 
to $\Sigma$. It is an {\em inner} one otherwise. 

The boundary $\partial I_{\mathcal R}$ of $I_{\mathcal R}$
is the set of boundary points in the totally geodesic 
hypersurface containing $I_{\mathcal R}$. 

A {\em pinched simple closed curve} is a simple curve
pinched at most three points or pinched at a connected arc. 
The boundary of the $I$-part is 
a disjoint union of pinched simple curves.


The following is a really important property
since this shows we can use crescents in general hyperbolic manifolds 
without worrying about whether the $I$-parts meet the boundary of
the ambient manifold. 

\begin{prop}\label{prop:disjcrescent} 
Let $\calR$ be a crescent in a $2$-convex ambient general hyperbolic 
manifold $M$. Then $\calR$ is disjoint from $\partial \tilde M$.
In fact $\calR$ is uniformly bounded away from 
$\partial \tilde M$.
\end{prop}
\begin{proof}
Suppose that $\calR$ meets $\partial \tilde M$. 
Since the closure of $\alpha_\calR$ being a subset of $\tilde \Sigma$ 
is disjoint from $\partial \tilde M$, 
it follows that $I_{\calR}$ meets $\partial \tilde M$ 
in its interior points and away from the boundary points 
in the ambient totally geodesic subsurface $P$ in $\tilde M$. 

We find the extreme point of $I_{\calR} \cap \partial M$ 
and find the supporting line. This point has a local half-space 
in $\calR$. By tilting the $I$-part a bit by the supporting line, 
we find a local half-space in $M$ and in it a local totally geodesic 
hypersurface meeting $\partial M$ at a point. This contradicts 
$2$-convexity. 
\end{proof}

\begin{figure}[ht]
\centerline{\epsfxsize=2.5in \epsfbox{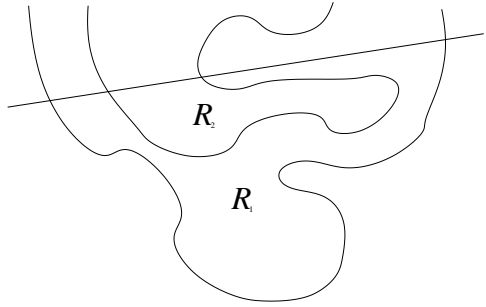}}
\caption{\label{fig:crescent} 
A $2$-dimensional section of crescents 
$\calR_1$ and $\calR_2$ where $\calR_2$ is $1$-nested 
and is innermost.}

\end{figure}
\typeout{<>}


The following shows the closedness of set of points of 
$\tilde M - \tilde \Sigma$ in crescents. 

\begin{prop}\label{prop:closedness} 
Let ${\calR}_i$ be a sequence of crescents. 
Suppose that $x$ is a point of $\tilde M - \tilde \Sigma$ 
which is a limit of a sequence of points in the union of 
${\calR}_i$. Then $x$ is contained in a crescent. 
\end{prop} 
\begin{proof} 
Let $x_i \in \bigcup_{j} {\calR}_j$ be a sequence 
converging to $x$. We may assume that $x$ is not an 
element of any ${\calR}_j$. 

Using geometric convergence, there exists a totally geodesic hypersurface
$P$ through $x$ and a geometric limit of a sequence of 
$I_{\calR_j}$ converging to a subset $D$ in $\tilde M$. 

Then $D$ is separating in $\tilde M$.
If not, there exists a simple closed curve $\gamma$ in $\tilde M$
meeting $D$ only once. This means $I_{\calR_j}$ for a sufficiently 
large $j$ meets $\gamma$ only once as well.

Now, $\tilde M - D$ may have more than one components. 
Since $x \in D$, we take a component whose closure $L$ 
contains $x$. $L \cap D$ is also separating.

\end{proof}


We may assume that the holonomy group of $\Sigma$ 
does not consist of parabolic or elliptic or identity 
elements only. 

A {\em size} of a crescent is the supremum of 
the distances $d(x, \alpha_{{\mathcal R}})$ for
$x \in I_{{\mathcal R}}$. We show that this is globally 
bounded by a constant depending only on $\Sigma$. 

\begin{prop}\label{prop:cresupper}
Let $M$ and $\Sigma$ be as above. 
Then there is an upper bound to the size of a crescent. 
\end{prop}
\begin{proof}
If not, using deck transformations acting on $\tilde \Sigma$,
we obtain a sequence of bigger and bigger 
compact crescents where the corresponding sequence of 
the $I$-parts leave any compact subset of $\tilde M$
and the corresponding sequence of $\alpha$-parts meets 
a fixed compact subset of $\tilde M$. 
Therefore, we form a subsequence of
the developing images of the $I$-parts converging to a point
of the sphere at infinity of $\hyps$. 

Let ${\mathcal R}_i$ be the corresponding crescents.
Then $\alpha_{{\mathcal R}_i}$ is a subsurface with boundary 
in the $I$-parts, and 
$\alpha_{{\mathcal R}_i}$ is getting larger and larger. 
 
Let $c$ be a closed curve in $\Sigma$ with nonidentity holonomy.
Let $\tilde c$ be a component of its inverse image in $\tilde \Sigma$.
Since $\tilde c$ must escape any compact subset of 
$\tilde \Sigma$, $\tilde c$ escape $\alpha_{{\mathcal R}_i}$. 
Thus, $\tilde c$ must meet all $I_{{\mathcal R}_i}$ for $i$ sufficiently large.
Since the developing image of $\tilde c$ has two well-defined endpoints, 
this means that the limit of the sequence of $I$-parts must 
contain at least two points, a contradiction. 
\end{proof} 

Given $\Sigma$, there is an upper bound to 
the folding-number of all crescents associated with $\Sigma$: 
This follows since $\tilde \Sigma$ is locally finite, 
and the increasing sequence of folding numbers implies that 
the sequence consists of crescents getting bigger and bigger. 

We call the maximum the {\em highest folding number} of $\Sigma$. 
We perturb $\Sigma$ to minimize the highest folding number 
which can change only by $\pm 1$ under perturbations.
After this, the folding number is constant under 
small perturbations. If there are no crescents, then 
the {\em folding number} of $\Sigma$ is defined to be $-1$.

Also, the union of all crescents for $\tilde \Sigma$ is 
in a uniformly bounded neighborhood of $\tilde \Sigma$
with the bound depending only on $\Sigma$.


We say that a $0$-folded crescent $\calR$ is 
a {\em highest-level} crescent 
if it is an innermost crescent of 
an $n$-folded crescent ${\calR}'$ where
$n$ is the highest-folding number of $\tilde \Sigma$ and $\calR$ is 
one that achieves the highest-level.

Suppose that $\calR$ is a compact highest-level crescent. 
Let $A_1, \dots, A_n$ be components of $I^O_{\calR}$. 
Recall that $I^O_\calR$ lies in a totally geodesic hypersurface.
The outermost pinched simple closed curve $\alpha_i$ 
in the boundary of $A_i$ has a trivial holonomy.
Since $\calR$ is of highest-level, 
$\alpha_i$ is an innermost curve itself 
or bounds innermost curves
in $\tilde \Sigma \cap \partial I^O_\calR$.
If each $\alpha_i$ is as in the former case,
then $\calR$ is said to be {\em innermost type crescent}, 
which is homeomorphic to a $3$-ball. 

%



We can classify the points of $\tilde \Sigma \cap I_{\mathcal R}$: 
A point of it is an {\em outer-skin point} if 
it is a nonboundary point and  
has a neighborhood in $\tilde \Sigma$ outside
${\mathcal R}^o \cup \alpha_{\mathcal R}$; 
a point is an {\em inner-skin point} if it is a nonboundary point and 
has a neighborhood in $\tilde \Sigma$ contained in 
${\mathcal R}$. A nonboundary point is either an outer-point
or an inner point or can be both.

The following classifies the set of outer-skin points. 
(A similar result holds for the set of inner points
except for (d).)

\begin{prop}\label{prop:sigmaI}
For a highest-level crescent $\calR$, the intersection points of 
$I^o_{\calR}$ and $\tilde \Sigma$ are either outer-skin points 
or inner-skin points. 
The set of outer-skin points of $I_{\mathcal R}$ for 
a highest-level crescent $\mathcal R$ is one of the following: 
\begin{itemize}
\item[(a)] a union of at most three isolated points. 
\item[(b)] a union of at most one point and a segment 
or a segment with some endpoints removed. 
\item[(c)] a union of two segments with a common endpoint
with some of the other endpoints removed. 
\item[(d)] a triangle with a boundary segment or 
two removed.
\end{itemize}
The same statement are true for inner-skin points.
\end{prop}
\begin{proof} 
This follows from the general position of vertices of $\tilde \Sigma$. 
\end{proof}

\begin{defn}\label{defn:I} 
Given a crescent $\calR$, we define $I^O_{\calR}$ to be the $I_{\calR}$
removed with the pinched points, boundary points, and the segments 
and triangles as above. (We don't remove the isolated points.)
\end{defn} 


\begin{prop}\label{prop:avoidg}
Suppose that $c$ is a straight closed geodesic in $M$ 
not meeting $\Sigma$. 
Let $\calR$ be a highest-level crescent.
Then 
\begin{itemize}
\item the inverse image $\tilde c$ of $c$ in $\tilde M$ 
does not meet $\calR$ in its interior 
and the $\alpha$-parts. 
\item $\tilde c$ could meet $\calR$ in its $I$-part
tangentially and hence be contained in the $I$-part.
In this case, $\calR$ is not compact.
\item If $\tilde c$ is a geodesic in $\tilde M$ 
eventually leaving all compact subsets, 
then the above two statements hold as well.
In particular if $\tilde M$ is a special hyperbolic 
manifold and $\tilde c$ ends in the limit set of 
the holonomy group associated with $\tilde M$. 
\end{itemize}
\end{prop}
\begin{proof}
If a portion of $\tilde c$ meets the interior of $\calR$, 
$\tilde c \cap \calR$ is a connected arc, say $l$
since $\calR$ is a closure of a component cut out by 
a totally geodesic hyperplane in $\tilde M - \tilde \Sigma$ --(*).

Since $\tilde c$ is disjoint from $\tilde \Sigma$, 
both endpoints of $l$ must be in $I^O_{\calR}$ or 
in $\SI^2_\infty \cap \ovl{I_{\calR}}$ for 
the closure $\ovl{I_{\calR}}$ of $I_{\calR}$ in 
the compactified $\hyps \cup \SI^2_\infty$. 
If at most one point of $l$ is in $I_{\calR}$, then 
$l$ is transversal to $I_{\calR}$ and the other endpoints 
$l$ must lie in $\alpha_{\calR}$ by (*). This is absurd. 

If at least two points of $l$ are in $I_{\calR}$, then 
$l$ is a subset of $I_{\calR}$. Since $\tilde c$ is disjoint from
$\tilde \Sigma$, $\tilde c$ is a subset of $I_{\calR}$, 
and $\calR$ is not compact. 

The only remaining possibility for $l$ is 
that $l$ ends in $\SI^2$ and $l^o$ is a subset of $\calR^o$. 
Hence $\tilde c = l$ and 
$\tilde c$ lies on the totally geodesic hypersurface containing 
$I_\calR$ since otherwise $l$ has to end at a point of $\alpha_\calR$. 

The domain $I_{\calR}$ in $P$ is disjoint from $\tilde c$ also on $P$.
Its boundary is a union of pinched arc in $\tilde \Sigma$. 
Let $\alpha$ be an open arc in $\hyps$ closest to $\tilde c$ on $P$. 
If an one-sided neighborhood of $\alpha$ goes into $\calR$,
Then $N(\alpha) - \alpha$ for a neighborhood $N(\alpha)$ 
of $\alpha$ lies outside $\calR$ since otherwise this leads to
a higher folding number.
Thus the vertices of $\alpha$ as an arc are vertices of 
$\tilde \Sigma$. Thus, $\alpha$ can have at most three vertices
as $\alpha$ is on a plane, and there must be an infinite 
geodesic edge of $\tilde \Sigma$. The later is impossible
since $\tilde \Sigma$ was finitely triangulated. 
Therefore, we obtained a contradiction, and $\tilde c$ can't 
be in the interior of $\calR$. This proves (i) and (ii). 

(iii) follows similarly.

\end{proof}


\begin{figure}[ht]
\centerline{\epsfxsize=3.8in \epsfbox{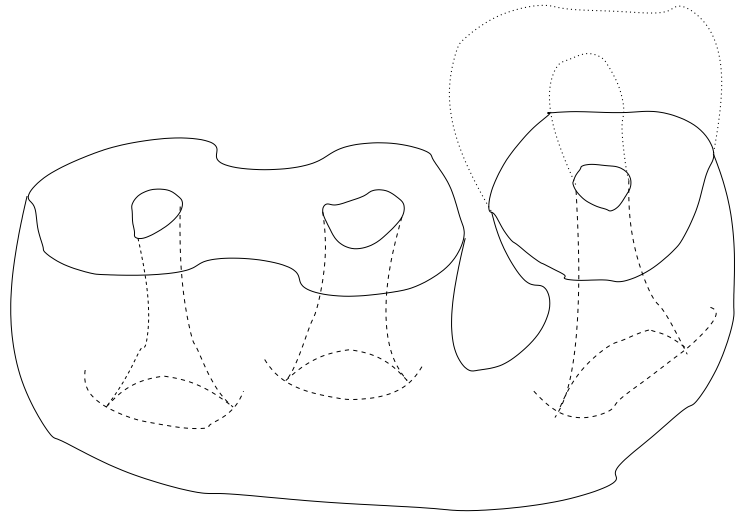}}
\caption{\label{fig:hlcrescent} The dashed arc indicates 
the tube from the bottom and the dotted arcs indicated 
the disks to be attached to the $\alpha$-parts}
\end{figure}
\typeout{<>}

\begin{prop}\label{prop:alphapt} 
Let $\calR$ be a highest-level crescent. Then 
there exists an innermost crescent 
$\calR'$ containing $\calR$ so that 
\begin{itemize}
\item $I_{\calR'}$ is connected 
and has no pinched points or a disconnecting union of 
outer-skin edges, i.e., $I^O_{\calR'}$ is connected.
\item The closure of $\alpha_{\calR'}$ is homeomorphic 
to $I_{\calR}$. 
\item $\calR'$ is homeomorphic to 
$I_{\calR'} \times [0,1]$.
\end{itemize}
\end{prop}
\begin{proof}
We assume the hypothesis \ref{hyp:incompressible}.

By Proposition \ref{prop:avoidg}, 
the interior of $\calR$ is disjoint from any 
lifts of $\cl_1, \dots, \cl_n$. 

Suppose that $\calR$ is compact to begin with.
Then $\calR$ is disjoint from the lifts of $\cl_1, \dots, \cl_n$
since $\tilde \Sigma$ is disjoint from these. 
We let $\calS$ be the highest-level crescent obtained 
from $\calR$ by cutting through the $I$-part hypersurface $P$
and taking the closure of 
a component of $\calR - \tilde \Sigma - P$ if necessary.

We introduce a height function $h$ on $\calS$ defined by 
introducing a parameter of hyperbolic hypersurfaces perpendicular 
to a common geodesic passing through $I_{\calS}$ 
in the perpendicular manner. (It will not matter which 
parameter we choose).
We may assume that $h$ is Morse in the combinatorial sense. 

If $I_\calS$ does not meet any lifts of $\cl_1, \dots, \cl_n$, 
let $N_\eps(\calS)$ be the neighborhood of $\calS$ in
the closure of the component of $\tilde M - \tilde \Sigma$ 
containing the interior of $\calS$. 

We let $N_{\eps}(\alpha_\calS)$ to be the intersection of 
$\tilde \Sigma$ with $N_{\eps}(\calS)$.
$N_\eps(\calS)$ can be chosen so that $N_\eps(\alpha_\calS)$ 
in $\tilde \Sigma$ becomes an open surface compactifying 
to a surface. There exists a part $I$ in of the boundary 
which is a complement in the boundary of $N_\eps(\calS)$ 
of $N_\eps(\alpha_\calS)$. We assume that $I$ lies on
a properly imbedded surface $P'$ perturbed from 
the $I$-part hypersurface $P$ of $\calS$. 

First, we show that $N_\eps(\calS)$ is a compression body 
with $N_\eps(\alpha_\calS)$ as the compressible surface in 
the boundary:

Topologically, $N_\eps(\alpha_\calS)$ is homeomorphic 
to a surface possibly with $1$-handles attached
from $\alpha_{\calS}$ and $I$ is obtained from $I_\calS$ by 
removing $1$-handles corresponding to the pinched points.

We may extend $h$ to an $\eps$-neighborhood of 
$\calS$. This may introduce only saddle type singularity 
in $N_\eps(\alpha_\calS)$. 
We let $I$ to be in the zero level of $h$.
This process only introduces a saddle type singularity in 
$N_\eps(\alpha_\calS)$.

If there is a critical point of $h$ with locally 
positive type, then we see that in fact 
there exists a crescent of higher level near the critical 
point. This is absurd. 

If there are no critical point of positive type, 
then $\pi_1(N_{\eps}(\alpha_{\calS})) \ra \pi_1(N_\eps({\calS}))$
is surjective as shown by Freedman-McMullen \cite{FM}. 

There exists a compression body in $N_{\eps}(\calS)$ 
with a boundary $N_{\eps}(\alpha_{\calS}) \cup S'$ 
for an incompressible surface $S'$ in the interior of $\calS$.
Since all closed path in $N_{\eps}(\calS)$ is 
homotopic to one in $N_{\eps}(\alpha_\calS)$, it follows that 
$S'$ is parallel to $I$. Hence, $N_{\eps}(\calS)$ is 
a compression body homeomorphic to $I$ times an interval 
and $1$-handles attached at disks disjoint from $I$. 
( $\calS$ is essentially obtained by pinching some points of $I$
together and pushing down a bit.) 

Next, we reduce the number of components of $I$:

Suppose now that $I$ is not connected. 
This means that there are $1$-handles attached to $I$ times the intervals 
joining the components. 
Then $N_{\eps}(\calS)$ has a compressing disk $D$ 
for $N_{\eps}(\alpha_\calS)$.
dual to the $1$-handles.
Since $\partial D$ bounds a disk $D'$ in $\tilde \Sigma$ by Dehn's lemma, 
$\partial D$ is separating in $\tilde \Sigma$.

Consequently also, $N_\eps(\alpha_{\calS})$ is a planar surface. 

The irreducibility of $\tilde M$ tells us that 
$D$ and $D'$ bound a $3$-ball $B$ in the closure of a component of 
$\tilde M - \tilde \Sigma$. Then $B$ contains at least one 
component of $I$. 

By taking a maximal family of compressing disks dual to the 
$1$-handles and regarding the components of the complements as 
vertices, we see that the $1$-handles do not form a cycle. 
Therefore, we choose the compressing disk $D$ to be 
the one such that $D$ and corresponding $D'$ bounds 
a $3$-ball $B$ containing a unique component of $I$. 

The ball $B$ contains the component of $\calR - D$ and 
a component of $I$, say $I'$.
Then $I'$ has a boundary 
component $\gamma$ in $D'$ and possibly other boundary components.
$\gamma$ bounds a disk $D''$ in 
$\tilde \Sigma - N_\eps(\alpha_\calS)$
which is in $D'$. The other boundary components
$\gamma_1, \dots, \gamma_m$ of $I'$ bound 
disjoint disks $D_1, \dots, D_m$ in the totally geodesic 
hypersurface $P'$ containing 
$I_{\calS}$. These $\gamma_i$s are innermost closed curves in 
$P' \cap \tilde \Sigma$. $I'$ union with these disks $D_i$s
and $D''$ bound a $3$-ball $B''$. $B'' \cap \tilde \Sigma$ 
is a union of $D''$ and a surface $\Sigma''$ with 
boundary equal to the union of the other closed boundary curves of $I'$.
Since $\Sigma'' \subset D'$, 
we have that $\Sigma''$ is a planar surface. 
If a component of $\Sigma''$ contains more that two of 
$\gamma_i$s, then we can find a closed curve in $D'$ meeting 
the two of $\gamma_i$s exactly once 
since $\alpha_{\calS}$ is connected. 
This is absurd. Therefore each component of $\Sigma''$ is a disk. 

Therefore, 
\[I \cup D'' \cup D''_1 \cup \cdots D''_m\] 
bounds a $3$-ball $B'''$. Taking a union of $N_\eps(\calS)$ with 
$B'''$, we obtain $N_\eps(\calT)$ for a crescent 
$\calT$ with $I_\calT$ in $I_\calS$ and one less component of 
the surface in $\partial N_\eps(\calT)$ corresponding to $I$.

By induction, we obtain a crescent $\calR''$ with 
$I_{\calR''}$ in $I_\calR$ and the surface $I'$ corresponding to 
$I$ connected. 
Since $\calR''$ is homeomorphic to a compression body,
$\calR''$ is homeomorphic to 
$I$ times an interval since there are no compressing disks. 

If there are any pinched points in $I_{\calR''}$ or 
disconnecting outer-skin edges, then $I'$ would be disconnected.
$\alpha_{\calR''}$ has a closure that is a surface since 
there are no pinching points.

Since $\calR''$ is an $I$-bundle, it follows that 
the closure of $\alpha_{\calR''}$ and $I_{\calR''}$ are 
homeomorphic surfaces. 

If there were more than one cut-off crescents $\calS$, then 
we obtain $\calR''$ for each $\calS$. 
There must be two cut-off crescents 
$\calS$ and $\calS'$ adjacent 
from opposite sides of some of the components of $I_\calS$.
Since the corresponding $\calR''$ and $\calR'''$  
does not have any pinched points or separating outer-skin edges, 
the unique components of $I_{\calR''}$ and $I_{\calR'''}$ either agree or 
are disjoint from each other. 
$\calR''$ and $\calR'''$ cannot be adjacent from opposite side
since we can form a compact component of 
$\tilde \Sigma$ otherwise. It follows that one of 
$\calR''$ and $\calR'''$ is a subset of the other. 
By induction, we see that the conclusions of the proposition
holds if $\calR$ is compact. 

If $\calR$ is noncompact, we follow as before but
we choose $N_{\eps}(\calR)$ similarly and
push down near infinity. Note here that only one 
component of $I$ maybe noncompact since the boundary of 
a compressing disk must bound a compact disk in $\tilde \Sigma$.
\end{proof}

\begin{defn} 
We say that $\calS$ as obtained from $\calR$ in 
the above lemma is obtained by cutting along 
the $I$-part hypersurface of $\calR$.  
$\calS$ may not be a unique one so obtained.
\end{defn}

\begin{cor}\label{cor:avoidg} 
Let $\calR$ be a secondary highest-level crescent.
Then the statements of Proposition \ref{prop:avoidg} 
hold for $\calR$ as well.
\end{cor}
\begin{proof}
The proof is exactly the same as that of Proposition \ref{prop:avoidg}.
\end{proof}

If a $0$-folded crescent $\calS$ contains $\calR$ so that $I^O_\calS$ is 
connected and is included in $I^O_\calR$, then 
we say that $\calS$ is a {\em highest-level} 
crescent as well. 
(Actually, it may not be highest-level since $\calS$ 
may not necessarily be contained in an $n$-folded crescent
but only a part of it.) 
More precisely, it is a {\em secondary highest-level crescent}.

Note that an outer highest-level crescent comes from 
an inner crescent if the highest folding number is odd
and comes from an outer crescent if the number is even.
The converse holds for an inner highest-level crescent. 


A secondary highest-level crescent exists for any 
highest-level crescent by Proposition \ref{prop:alphapt}:

\begin{cor}\label{prop:topproperty} 
Let $\calR$ be a highest-level crescent. 
If $I^O_{\calR}$ is not connected, 
then there exists a $0$-folded crescent $\calS$ containing $\calR$ 
so that $I^0_{\calS}$ is connected and 
is a component of $I^O_{\calR}$. 
\end{cor}

Also, a secondary highest-level crescent
has the outer-skin points with the same properties 
as those of a highest-level crescent.
(See Proposition \ref{prop:sigmaI}).

By taking a nearby crescent inside, we see that 
a highest-level crescent could be {\em generically chosen}
so that the crescent is compact,
the $I$-part and the $\alpha$-part are surfaces,
and $I^O$-part is truly the interior of the $I$-part. 

\begin{cor}\label{cor:highestcrestop} 
Let $\calR$ be the compact secondary highest-level 
outer (resp. inner) crescent that is generically chosen.
Then $\calR$ is homeomorphic to the closure of $\alpha_{\calR}$ 
times $I$, and $I^O_{\calR}$ is isotopic to $\alpha_{\calR}$ 
by an isotopy inside $\calR$ fixing the boundary of $I_{\calR}$. 
\end{cor}



We say two crescents $\calR$ and $\calS$ {\em face each other} 
if $I_{\calR}$ and $I_{\calS}$ agree with each other in some 
$2$-dimensional part and have disjoint one-sided neighborhoods.

\begin{prop}\label{prop:nofacing} 
If two highest-level outer 
{\rm (}resp. inner{\rm )} crescents $\calR$ and $\calS$ face 
each other, then there exists a (secondary) highest-level
outer {\rm (}resp. inner{\rm )} crescent $\calT$ 
with connected $I^O_{\calT}$ containing both. 
\end{prop} 
\begin{proof} 
We may assume without loss of generality that
$I^O_{\calR}$ and $I^O_{\calS}$ are both connected. 
If not, we replace $\calR$ and $\calS$ by ones 
with connected $I^O$-parts. The replacements
still face each other or one becomes a subset of 
another. In the second case, we are done. 

Since $I^O_\calR$ and $I^O_\calS$ meet in open subsets, 
either they are identical or the boundary $\partial L$ of 
their intersection $L$ lies $I^O_\calR$.
$\partial L$ is a subset of $\tilde \Sigma$ and 
is a $1$-complex of pinched arcs.
$\partial L$ is a set of outer-skin points of $\calS$ 
since a neighborhood of $\partial L$ in $\tilde \Sigma$ must 
be above $\calS$. However, then $\partial L$ must 
be disjoint from $I^O_\calS$. 
Therefore, $I^O_\calS = I^O_\calR$. 
This means that $\calS \cup \calR$ is bounded by 
a component subsurface of $\tilde \Sigma$.
This is absurd.

\end{proof}


\begin{defn}\label{defn:transvint} 
Two secondary highest-level outer (resp. inner) crescents 
$\calR$ and $\calS$ are said to meet 
{\em transversally} if $I_{\calR}$ and $I_{\calS}$ meet in 
a union of disjoint geodesic segment $J$, $J \ne \emp$, mapping into 
a common geodesic in $\hyps$, in a transversal manner such that 
\begin{itemize}
\item The the closure $\nu_{\calR}$ of the union of 
the components of $I_{\calR} - J$ 
in one-side is a subset of $\calS$ 
and the closure $\nu_{\calS}$ of the union of those of 
$I_{\calS} - J$ is a subset of $\calR$. 
\item The intersection ${\calR} \cap {\calS}$ is the closure of 
$\calS - \nu_{\calR}$ and conversely the closure of 
$\calR - \nu_{\calS}$. 
\item The intersection $\alpha_{\calR} \cap \alpha_{\calS}$ is 
a union of components of $\alpha_{\calR} - \nu_{\calS}$ in 
one-side of $\nu_\calS$
and, conversely, is a union of components of 
$\alpha_{\calS} - \nu_{\calR}$ in one side of $\nu_\calR$.
\item $\alpha_{\calR} \cup \alpha_{\calS}$ is an open surface in 
$\tilde \Sigma$. 
\end{itemize}
\end{defn}


\begin{prop}\label{prop:transvint} 
Given two secondary highest-level 
outer (resp. inner) crescents $\calR$ and $\calS$, 
there are the following mutually exclusive possibilities: 
\begin{itemize}
\item $\calR$ and $\calS$ do not meet in $\tilde M -\tilde \Sigma$. 
\item $\calR \subset \calS$ or $\calS \subset \calR$. 
\item $\calR$ and $\calS$ meet transversally. 
\end{itemize}
\end{prop}
\begin{proof} 
The reasoning is exactly the same as \cite{cdcr1} and \cite{psconv} 
in dimension two or three.
\end{proof}

\section{Crescents and isotopy}


In section 3, the crescent-isotopy theory is presented:
We form the union of secondary highest-level crescents
and can isotopy the union of 
their $\alpha$-parts to the complement $I$ in 
the boundary of their union. 
However, there might be some parts of $\tilde \Sigma$ 
meeting $I$ tangentially from above. We need 
to push these parts upward first
using truncations before the actual crescent-isotopy.

We first discuss how small isotopy can affect the 
set of all crescents of $\tilde \Sigma$.
For a crescent in the outer-direction, the move 
in the outer-direction ``preserves'' the crescent 
by moving the $\alpha$-part only.
For a crescent in the inner-direction, the move 
in the inner-direction also preserves. 

We introduce truncation move: we truncate a convex 
vertex by a local totally geodesic hypersurface and perturb
the result into general position by pushing the vertices. 
We show that the crescents in these cases are preserved or 
we can modify the crescents by moving the $I$-parts. 
These moves are designed to eliminate 
some points so that the result of the isotopy are 
imbeddings. 

We start from the outer-direction secondary highest-level crescents.
We take a union of an overlapping equivalence class of them 
and show that the union of 
$\alpha$-parts in $\tilde \Sigma$ can be isotopied to 
the complement in the boundary of the union. 
By above truncation moves, we show that the result is 
an imbedding. Finally we do this for all the equivalence 
classes and we can isotopy $\Sigma$ itself. 

The inner-direction crescent moves are entirely the same.


\begin{lem}\label{lem:crescentiso}
Suppose that $\tilde \Sigma$ has been isotopied in the outward direction 
by a sufficiently small amount
and $\calR$ is an outer crescent. Then there exists 
a crescent $\calR'$ sharing the $I$-part hypersurface with 
$\calR$ and differs from $\calR$ by isotopying the 
$\alpha$-part only. 
Conversely, if $\tilde \Sigma$ has been isotopied in the inward direction 
and $\calR$ is an inner crescent, the same can be said. 
\end{lem}
\begin{proof} 
Straightforward.
\end{proof}

We say that $\calR'$ is {\em isotopied from $\calR$ with the 
$I$-part preserved}.



We may ``truncate" $\Sigma$ at convex vertices and 
$\tilde \Sigma$ correspondingly and perturb: 
Let $v$ be a convex vertex and $H$ a local half-open ball 
at $v$ with the side $F$ passing through $v$. 
We may move $F$ inside by a very small amount 
and then truncate $\tilde \Sigma$ and correspondingly 
for all vertices equivalent to $v$.
Then we introduce some equivariant triangulation of 
the trace $T$ of the truncation and the truncated $\tilde \Sigma$
without introducing vertices in the interior of $T$. 
We push the concave vertices of $T$ downward by small
amounts along the corresponding edges of $\tilde \Sigma$ 
and then move the vertices inward
to make the truncated surface to be in general position.
The three steps together are called the {\em small truncation move}.


We denote by $\Sigma^\eps$ the perturbed $\Sigma$ where 
the traces are less than an $\eps$-distance away from the 
respective convex vertices. 
We assume that during the perturbations $\Sigma^\eps$ is 
isotopied from $\Sigma$ and the convexity of the dihedral
angles do not change under the isotopy. 
Thus, if an edge or a vertex is convex after being born, it will
continue to be so as $t \ra 0$ and as $t$ grows from $0$.

We may view the truncation move by considering $v$ 
to be a vertex of some multiplicity and vertices diverge as 
the side $F$ moves away from $v$. That is, we
see this as births of many vertices from convex vertices.


We may also assume that the convex vertex move is 
equivariant on $\tilde \Sigma$, i.e., the isotopy 
is equivariant. 

An {\em isotopy} of a crescent as we deform $\Sigma$ 
is a one-parameter family of crescents $\calR_t$ with 
$\alpha$-parts in $\Sigma$.
We say that a crescent {\em bursts} if fixing the totally geodesic 
hypersurface containing the $I$-part of it and isotopying
the $\alpha$-parts in the isotopied $\Sigma$ 
cannot produce a crescent isotopied from the original one. 


Such an event happens when a parameter of an edge of $\tilde \Sigma$ 
or a parameter of a vertex of $\tilde \Sigma$ 
go belows the fixed totally geodesic hypersurface from the 
point of view of the crescent. 
Of course, a vertex could be a multivertex and 
all of the new vertices go down. 
The edge should be one the face that meets the $I$-parts of
the crescents and the vertex on the edge that meets the 
$I$-part of the crescent.
The event could happen simultaneously but 
the generic nature of the move shows that 
at most four vertex submersions, at most three edge submersions, 
and at most two vertices and one edge submersions can 
happen simultaneously. 
Moreover, at the event, 
the vertex and the edge must be in the $I$-part of 
the crescent and the triangles of $\Sigma$ must be placed 
in certain way in order that the bursting to take place.

If the bursting happens immediately for a crescent $\calR$, 
then the convex vertex must be in the boundary of 
$I_{\calR} \cap (\tilde M - \tilde \Sigma)$ in 
the totally geodesic hypersurface $P$ containing $I_{\calR}$.
Otherwise a small perturbation of $\tilde \Sigma$
gives us an isotopy of $\calR$ preserving the $I$-part hypersurface. 
Therefore, the bursting does not take place.

\begin{prop}\label{prop:perturb}
Under a small truncation move in the outer direction, 
we can isotopy
\begin{itemize}
\item[(i)] each outer crescent into itself by moving 
the $\alpha$-part in the outer direction
and preserving the $I$-part hypersurface. 
\item[(ii)] each inner crescent into itself union 
the $\eps$-neighborhood of $\tilde \Sigma$ by moving 
the $I$-part hypersurface in the outer direction or 
preserving the $I$-part hypersurface.
\end{itemize} 
Under a small truncation move in the inner direction, 
we can deform 
\begin{itemize}
\item[(iii)] each inner crescent into itself by moving 
the $\alpha$-part in the inner direction and 
preserving the $I$-part hypersurface.
\item[(iv)] each outer crescent into itself union 
the $\eps$-neighborhood of $\tilde \Sigma$ by moving 
the $I$-part hypersurface in the inner direction
or preserving the $I$-part hypersurface. 
\end{itemize} 
All crescents of $\tilde \Sigma^\eps$
can be obtained in this way.
The highest folding number may decrease only
under a convex vertex move.
\end{prop}
\begin{proof} 
Essentially, the idea is that the move can only ``decrease" 
the associated crescents. 


Let $\calR$ be an outer crescent and $\tilde \Sigma$ moved 
in the outer direction. Lemma \ref{lem:crescentiso} implies (i).

Let $\calR$ be an inner crescent and $\tilde \Sigma$ be moved 
in the outer direction. Then again an isolated submerging 
vertex is a convex vertex.  In this case, we move the 
$I$-part inward so that the submerging vertex stay on
the boundary of the $I$-part. Other cases are treated similarly.
This proves (ii).

(iii) and (iv) correspond to (i) and (ii) respectively
if we change the orientation of $\Sigma$.

To show that all crescents of $\tilde \Sigma^\eps$
can be obtained in this way:
Given an outer crescent for $\tilde \Sigma^\eps$, 
we reverse the truncation move. 
If the $I$-part of a crescent avoids the 
traces of the truncation moves, then we simply
isotopy the $\alpha$-parts only. 

The trace surface has only concave vertices and 
s-vertices. Let $P'$ be a local totally geodesic hypersurface 
truncating the stellar neighborhood of a convex vertex 
$v$ of $\tilde \Sigma$ at some small distance from $v$ 
but large compare to our isotopy move distance.
Suppose that $v$ were involved in the convex truncation move.  

Suppose that the $I$-part of a crescent $\calR$ for $\tilde \Sigma$ 
meets what are outside the part truncated by $P'$. 
Assuming that our isotopy was very small, 
if the $I$-part meets one of the trace surface, 
the $I$-part meets $P'$ since $P'$ is separating.
$P'$ intersected with the closure of the exterior of $\tilde \Sigma^\eps$ 
is a polygonal disk $D^\eps$.
Then $D^\eps$ intersected with the $I$-part is a disjoint union
of segments. Letting $D$ be the intersection of $P'$ 
with the exterior of $\tilde \Sigma$, we see that 
$D$ intersected with $P'$ extend the segments of 
$D^\eps$ intersected with $P'$. 
Thus, it is clear that the $I$-part extends into the polyhedrons 
bounded by $P'$ and the stellar neighborhood of $\tilde \Sigma$.
Since all vertex submersions of $\tilde \Sigma^\eps$ 
can happen by vertices near the convex vertices of
$\tilde \Sigma$ masked off by 
totally geodesic hypersurfaces such as $P'$,
we obtain a crescent $\calR'$ for $\tilde \Sigma$
preserving the $I$-part hypersurface of $\calR$. 
Therefore, $\calR$ were obtained from $\calR'$
by the convex truncation isotopy preserving the $I$-part.

Therefore a crescent for $\tilde \Sigma^\eps$ is one 
we obtained by the process in (i). 

Let $\calR$ be an inner crescent for $\tilde \Sigma^\eps$.
Then since the vertices moved outward with respect to 
$\tilde \Sigma$, they move inward when we reverse the process and 
we see that $\calR$ is isotopied to a crescent for
$\tilde \Sigma$ by Lemma \ref{lem:crescentiso}
by preserving the $I$-part hypersurface.

To show that the highest folding number can only decrease:
For a crescent $\calR$ to increase the folding number, 
a vertex must move into $I_{\calR}$ during the isotopy. 
We see that such a vertex must be a convex one. 
However, the convex vertices can only move in the 
direction away from the interior of $\calR$. (Even ones 
after the births obey this rule.)

\end{proof}


First, we suppose that there are highest-level crescents that are 
outer ones. We move $\tilde \Sigma$ in the outer direction first to 
eliminate the outer highest-level crescents.


As we did in \cite{cdcr1} and \cite{psconv}, 
we say that two highest-level crescents $\calR$ and $\calS$ 
are equivalent if there exists a sequence of 
transversally intersecting crescents from 
$\calR$ to $\calS$; that is, 
\[\calR = \calR_0, \calR_i \cap \calR_{i+1}^o \ne \emp, 
\calS = \calR_n \hbox{ for } i = 1, 2, \dots, n. \]

We define $\Lambda({\calR})$ to be the union of 
all highest-level crescents equivalent to 
the highest-level crescent $\calR$. 
As before, $\Lambda({\calR})$ and $\Lambda({\calS})$ 
do not meet in the interior or they are the same. 

We define $\partial_I \Lambda({\calR})$
to be the boundary of $\Lambda({\calR})$ removed 
with the closure of the union of the $\alpha$-parts of 
the crescents in it. 
Then $\partial_I \Lambda({\calR})$ is 
a convex surface. 

We define $\partial_\alpha \Lambda({\calR})$ 
as the union of the $\alpha$-parts of the crescents 
equivalent to $\calR$. 

Recall that a pleated surface is a surface where 
through each point passes a geodesic. 

\begin{lem}\label{lem:pIlambda} 
The set 
\[\partial_I \Lambda({\calR}) \cap \tilde M - \tilde \Sigma\]
is a properly imbedded pleated surface.
\end{lem}
\begin{proof}
For each point of $x$ belonging to the above set, 
$x$ is an element of the interior of $\tilde M$ 
by Proposition \ref{prop:disjcrescent}. 
Let $B(x)$ be a small convex open 
ball with center at $x$. Then the 
crescents equivalent to $R$ meet $B(x)$ 
in half-spaces. Therefore the complement of 
their union is a convex subset of $B(x)$
and $x$ is a boundary point. There is a supporting 
half-space $H$ in $x$ and $H$ belongs to 
$\Lambda(\calR)$. 

If there were no straight geodesic 
passing through $x$ in the boundary 
set $\partial_I \Lambda(\calR)$, then
there exists a totally geodesic disk $D$ in $B(x)$  
with $\partial D$ in $\Lambda(\calR)$ 
but interior points are not in it. 

Since $\partial D$ is in $\Lambda(\calR)$, 
each point of $\partial D$ is in some crescent. 
We can extend $D$ to a maximal totally geodesic hypersurface 
and we see that a portion of the hypersurface
bounds a crescent $\calT$ containing $D$ in its $I$-part 
and overlapping with the other crescents.
Thus $\calT$ is a subset of $\Lambda(\calR)$
and so is $D$. 

Therefore, $\partial_I \Lambda(\calR)- \tilde \Sigma$ 
is a pleated surface. 
\end{proof}


We may have some so-called outer-skin points of $\tilde \Sigma$ at
$\partial_I \Lambda({\calR})$, i.e, those points with neighborhoods 
in $\tilde \Sigma$ outside $\Lambda(\calR)$. 
We can classify outer-skin points.

\begin{prop}\label{prop:outerpts}
The set of outer-skin points on $\partial_I \Lambda(\calR)$ 
is a union of the following components: 
\begin{itemize}
\item isolated points, 
\item an arc passing through the pleating locus
with at least one vertex. 
\item isolated triangles, 
\item union of triangles meeting each other 
at vertices or edges. 
\end{itemize}
\end{prop}
\begin{proof} 
This essentially follows by Proposition \ref{prop:sigmaI}.
\end{proof}

There have to be convex vertices in all of the above 
cases of outer-skin points. 
We move a convex vertex by vertex in an equivariant 
manner. Therefore, we may move them inwardly by small truncation moves, 
i.e., along an inward normal to the local half-open balls. 
This will only decrease the set of the union of crescents. 
In every case, a new convex vertex that is an outer-skin point 
is obtained after the move.

In this way, we see that 
\begin{equation}\label{eqn:replace}
\partial_I \Lambda({\calR}) 
\cup (\tilde \Sigma - \partial_\alpha \Lambda({\calR})) 
\end{equation} 
is a properly imbedded pleated surface. 

We do this for $\Lambda({\calR})$ for each highest-level crescent 
$\calR$. The end result is a properly imbedded surface 
$\tilde \Sigma'$. The deck transformation group 
acts on $\tilde \Sigma'$ since it acts on the union of 
$\Lambda({\calR})$. Thus, we obtain a new closed surface 
$\Sigma'$. 

Since the union of $\Lambda({\calR})$ for every highest-level crescent 
$\calR$ is of bounded distance from $\tilde \Sigma$ 
by the boundedness and the fact that $M$ is locally-compact, 
$\Sigma'$ is a compact surface.

We show that $\Sigma$ and $\Sigma'$ are isotopic. 

Let $N$ be the $\eps$-neighborhood of $\tilde \Sigma'$ in 
the closure of the outer component of $\tilde M - \tilde \Sigma$. 
There exists a boundary component $\partial_1 N$ nearer 
to $\tilde \Sigma$ than the other boundary component. 
The closure of a component $K$ of $\tilde M - \tilde \Sigma - \tilde \Sigma'$
contains $\partial_1 N$. 
Then $K$ projects to a compact subset of $M$. 
We can find a finite collection of generic secondary highest-level 
compact crescents 
$\calR_1, \dots, \calR_n$ and whose images under 
$\Gamma$ form a locally finite cover of $K$. 

We label the crescents by $\calS_1, \calS_2, \dots$. 
We know that replacing the closure of the $\alpha$-part of 
$\calS_1$ by the $I$-part is an isotopy. 
After this move, $\calS_2, \calS_3, \dots$
become new generic highest-level crescents 
by Proposition \ref{prop:transvint} and appropriate 
truncations. 

We define $\partial_I (S_1 \cup S_2 \cup \dots)$ 
as the boundary of $S_1 \cup S_2 \cup \dots$ 
removed with the union of the $\alpha$-parts of $S_1, S_2, \dots$. 
Again, this is a convex imbedded surface. 
Therefore, replacing the union of 
the $\alpha$-parts of $\calS_1, \calS_2, \dots$ 
by $\partial_I (S_1 \cup S_2 \cup \dots)$ is 
an isotopy as above. 

We obtain $\Sigma_{\calR_1, \dots, \calR_n}$ as 
the image in $M$, which is isotopic to $\Sigma$.
If $\eps$ is sufficiently small, then 
we see easily that 
$\partial_I (S_1 \cup S_2 \cup \dots)$ is in $N$ 
can can be isotopied to $\tilde \Sigma'$ 
using rays perpendicular to $\tilde \Sigma'$. 
Thus, $\Sigma_{\calR_1, \dots, \calR_n}$
is isotopic to $\Sigma'$.
We showed that $\Sigma'$ is isotopic to $\Sigma$.

\section{Convex perturbations} \label{sec:convpert}

In this section, we discuss how a surface with a portion of 
itself concavely pleated by infinitely long geodesics and the remainder 
triangulated can be perturbed to a triangulated surface 
without introducing higher-level crescents. 
This is done by approximating the union of pleating geodesics
by train tracks and choosing normals in the concave direction 
and finitely many vertices at the normal and pushing 
the pleating geodesics to become geodesics broken at 
the vertices.


We will now perturb the isotopied $\Sigma$ into 
a triangulated surface not 
introducing any higher-level crescents. 

Suppose that $\Sigma$ is a closed imbedded surface in $M$.
$\Sigma$ is a {\em pleated-triangulated surface} 
if $\tilde \Sigma$ contains 
a closed $2$-dimensional domain divided into locally finite collection
of closed totally geodesic convex
domains meeting each other in geodesic segments
and through each point of the complement passes 
a {\em straight} geodesic.
We may also add finitely many straight geodesic segments 
in the surface ending at the domain.
The domain union with the segments 
is said to be the {\em triangulated part} of 
$\Sigma$. 
The boundary of the domain is a union of finitely pinched 
simple closed curves. 
The complement of the domain is an open surface, which is
said to be the {\em pleated part} where through each point 
passes a straight geodesic.
Then the pleated part has a locus where a unique 
straight geodesic passes through. This part is said to be the 
{\em pleating locus}. It is a closed subset of the complement and 
forms a lamination. 

For later purposes, we say that $\Sigma$ is 
{\em truly pleated-triangulated} if 
the triangulated part is a union of totally
geodesic domains that are polygons and geodesic 
segments ending in the domains.

Note that the triangulated parts and pleating parts are 
not uniquely determined. We simply choose what seems 
natural.

We assume that the pleated part is locally convex or 
locally concave. We choose a normal direction 
so that the surface is locally concave at the pleated part. 

If $\Sigma$ satisfies all of the above conditions, 
we say that $\Sigma$ is a {\em concave pleated-triangulated} 
surface. If we choose the opposite normal-direction, 
then $\Sigma$ is a {\em convex pleated, triangulated} 
surface.



First, we show that the surface obtained in the previous 
section is concave pleated-triangulated. 

\begin{prop}\label{prop:ptsurface}
Let $\Sigma$ be a triangulated surface incompressible 
in $M$ with a finite number of closed geodesics removed.
Let $n$ be the highest folding number of $\Sigma$. 
Suppose $n$ is achieved by outer crescents. 
Suppose that $\Sigma'$ is obtained from $\Sigma$
by a small truncation move for level $n$ outer crescents
and next by crescent move for level $n$ outer crescents. 
Then $\Sigma'$ is a concave pleated-triangulated 
surface in the outer direction.
Moreover, the statements are true if all ``outer" were 
replaced by ``inner".
\end{prop}
\begin{proof} 
The part $\partial_I \Lambda(\calR) - \tilde \Sigma$ 
for crescents $\calR$s are pleated by Lemma \ref{lem:pIlambda}. 
These sets for crescents $\calR$ are either identical or 
disjoint from each other as the sets of form 
$\Lambda(\calR)$ satisfy this property. 
The union of sets of form $\partial_I \Lambda(\calR)$ 
form the pleated part and the complement in 
$\Sigma'$ were in $\Sigma$ originally and 
they are union of totally geodesic 
$2$-dimensional convex domains. 
\end{proof} 





Next, we analyze the geometry of the pleating locus 
in the pleated part.

Two leaves are {\em converging} if one 
is asymptotic to the other one (see Section \ref{sec:metricspaces}
for definitions).
i.e., the distance function from one leaf to the other 
converges to zero and conversely.
By an end of a leaf of a lamination {\em wrapping around} 
a closed set, we mean that the leaf converges to a subset of 
the closed set in the direction of the end. 


\begin{lem}\label{lem:nowrapping} 
Suppose that $\Sigma$ is a closed concave pleated-triangulated 
surface with a triangulated part and pleated part assigned. 
Suppose that $l$ is a leaf.
Then 
\begin{itemize}
\item $l$ is either a leaf of a minimal geodesic lamination or 
a closed geodesic, or each end of $l$ wraps around  
a minimal geodesic lamination or a closed geodesic
or ends in the boundary of the triangulated part.
\item if $l$ is isolated from both sides, then 
$l$ must be a closed geodesic, and 
\item pleating leaves in a neighborhood of $l$ 
must diverge from $l$ eventually.
\end{itemize}
\end{lem} 
\begin{proof}
Since the pleated open surface carries an intrinsic metric 
which identifies it to a quotient of 
an open subset of the hyperbolic space, 
each geodesic in the pleating lamination will 
satisfy the above properties like the geodesic laminations 
on the closed hyperbolic surfaces.

If $l$ is isolated, then there is a definite positive 
angle between two totally geodesic hypersurfaces ending at $l$.
Suppose that $l$ is not a simple closed geodesic. 
Then this angled pair of the hypersurfaces continues to wrap around
infinitely in $M$ accumulating at a point of $M$ and 
the sum of the angles violates the imbeddedness of $\Sigma$.

If $l$ is not isolated but has converging 
nearby pleating leaves, the same reasoning will hold. 
Therefore, the second item holds. 
\end{proof}

\begin{prop}\label{prop:pleatinglocus} 
Suppose that $\Sigma$ is a closed concave pleated-triangulated 
surface. Let $\Lambda$ be the set of pleating locus of the pleated
part in $\Sigma$. Then $\Lambda$ decomposes into finitely many 
components $\Lambda_1, \dots, \Lambda_n$ so 
that each $\Lambda_i$ is one of the following:
\begin{itemize}
\item a finite union of finite-length pleating leaves homeomorphic to 
a compact set times a line. {\rm (} A discrete set times a line 
if $\Sigma$ is truly pleated-triangulated. {\rm )}
\item a simple closed geodesic.
\item a minimal geodesic lamination, which is a closed 
subset of the pleated part isolated away from the triangulated 
part. 
\end{itemize} 
Here each leaf is either bi-infinite or finite.
The union of bi-infinite leaves is a finite union of 
minimal geodesic laminations and is isolated way 
from the triangulated part and the union of 
finite-length leaves. 
\end{prop}
\begin{proof}
Let $l$ be an infinite leaf in the pleating locus. 
By Lemma \ref{lem:nowrapping}, $l$ is not 
isolated from both sides and the leaves in its neighborhood 
is diverging from $l$. 
If $l$ is not itself a leaf of a minimal lamination, 
then an end of $l$ must converge to a minimal lamination or 
a simple closed geodesic. This means that leaves in 
a neighborhood also converges to the same lamination
in one of the directions. However, this means that they 
also converges to $l$, a contradiction. 
Therefore, each leaf is a leaf of a minimal lamination or 
a closed geodesic or a finite length line. 

The union of all finite length lines in the pleating locus 
is a closed subset: Its complement in $\Lambda$ is 
a compact geodesic lamination in $\Sigma$.
If a sequence of a finite length leaves $l_i$ converges to 
an infinite length geodesic $l$, then $l_i$ gets 
arbitrarily close to a minimal lamination or
a closed geodesic. If $l_i$ gets into a sufficiently
thin neighborhood of one of these, then
an endpoint of $l_i$ must be in a sufficiently 
thin neighborhood of one of these by
the imbeddedness property of $l_i$, i.e., 
cannot turn sharply away and go out of the neighborhood.
As $l_i$ ends in the triangulated part, 
the distance from the triangulated part to 
one of these goes to zero. 
Since the domains in the triangulated part 
are in general position, the boundary of
the triangulated part cannot contain 
a closed geodesic or the straight geodesic lamination. 
This is a contradiction.

Looking at an $\eps$-thin neighborhood of $\Lambda$, 
we see that $\Lambda$ decomposes as described.
(See Casson-Bleiler \cite{Cass-Bl} for background informations).
\end{proof}

\begin{rem} If $\Sigma$ is properly pleated-triangulated, 
then there are only finitely many finite length 
pleating leaves since their endpoints are on 
the vertices of the triangles in the triangulated part. 
\end{rem}


We can still define crescents for $\tilde \Sigma'$ which is 
the inverse image of $\Sigma'$ obtained from crescent moves. 
The definitions are of course the same. 

\begin{prop}\label{prop:nocrescent} 
If $\Sigma'$ is obtained from $\Sigma$ by 
the highest-level outer crescent move for level $n$, 
then the union of the collection of crescents of $\tilde \Sigma$ 
contains the union of those of $\tilde \Sigma'$, and 
$\Sigma'$ has no outer crescent of level $n$ or higher.
The same statements hold if we replace the word ``outer'' by 
``inner''.
\end{prop}
\begin{proof}
The outer crescents for $\tilde \Sigma'$ 
can be extended to ones for $\tilde \Sigma$ since 
their $I$-part can be extended across. 

The inner crescents for $\tilde \Sigma'$ can be 
truncated to ones for $\tilde \Sigma$ by Lemma \ref{lem:crescentiso}
since the move from $\tilde \Sigma'$ to $\tilde \Sigma$ 
is inward and supported by the outer crescents of $\tilde \Sigma$.
Thus the first statement holds. 

If there were outer crescent $\calR$ of level $n$ or higher, 
then we can extend $I_\calR$
across the $I_\partial$-parts so that we can obtain 
a level-$n$ or higher-level crescent. This is absurd. 
If $\calR$ were inner, Lemma \ref{lem:crescentiso}
implies the result.
\end{proof}

There might still be level $n-1$ outer crescents for 
$\Sigma'$. We do small truncation moves for these. 


We will use the train tracks to prove the following theorem:

\begin{thm}\label{thm:convexperturb} 
Let $\Sigma$ be a closed concave pleated, triangulated surface
with outer {\rm (} resp. inner {\rm )} level strictly less than $n$.
Then one can find an imbedded isotopic 
triangulated surface $\Sigma'$ in 
any $\eps$-neighborhood of $\Sigma$ so that
the following hold. 
\begin{itemize}
\item[(i)] The union of the set of crescents for $\tilde \Sigma'$ 
is in the $\eps$-neighborhood of that of $\tilde \Sigma$
and vice-versa for a small $\eps$ 
if we choose $\tilde \Sigma'$ sufficiently close to $\tilde \Sigma$.  
\item[(ii)] Assuming that $\Sigma$ is truly pleated-triangulated,
the vertices of the triangulated part 
of $\Sigma'$ are deformed from the vertices that $\Sigma$ had
and the vertices formed from points of the pleated parts.
{\rm (} The pleated part vertices are strict s-vertices.{\rm )} 
\end{itemize}
\end{thm}


A train track is obtained by taking a thin neighborhood of 
the lamination. We can think of the train track as a union of 
segment times an interval, so-called branches, joined up at the end of 
each segment times the intervals so that the intervals 
stacks up and matches. A point times the interval 
is said to be a {\em tie} and a tie where 
more than one branches meet a {\em switch}.
One can collapse the interval direction to obtain 
a union of graphs and circles.


\begin{proof}
Here, we will be working in $M$ directly.

(ii) Let $l_1, \dots, l_k$ be the thin strips 
containing all the 
finite length open leaves. 
We find a thin totally geodesic hypersurface 
$P_i$ near $l_i$s normal to the normal vectors of $I_i$s. 
Then we cut off the neighborhood of $l_i$ in $\Sigma$ 
by $P_i$ and replace the lost part with the portion in
$P_i$. We triangulate the portion and perturb the vertex 
inward to obtain a surface whose inverse image 
in $\tilde M$ is in a general position. 
This introduces squares which are triangulated into
pairs of triangles. 

This forms a generalization of a small truncation move. 
We still call it a small truncation move. --(*)

If there are outer level $(n-1)$-crescents with 
outer-skin points, then we remove the outer-skin points by 
small truncation-moves at convex vertices. --(**)

We now add finite leaves of infinite length in 
the pleated part so that the components of the 
complement of the union of the pleating locus and 
these leaves are all open triangles.
This can be done even though the boundary of
the pleated part is not geodesic. 

By choosing sufficiently small $\eps$-thin-neighborhoods
of the union for $\eps>0$,  we obtain a train track. 
We assume that the circle component contains 
a unique closed geodesic. 


We first choose switches for the endpoint of the finite 
length geodesics in the squares and added infinite length
finite leaves. We choose switches for the rest. 
We label them $I_1, \dots, I_n$.
The train track collapses to a union of 
graphs with vertices corresponding to $I_1, \dots, I_n$ 
and closed geodesics. 
The complement is a disjoint union of open triangles.

By choosing $\eps >0$ sufficiently small, we can assume that 
the outer-normal vectors to totally geodesic hypersurfaces near
$I_i$ are $\delta$-close for a small $\delta > 0$
except the outer-normal vectors to the totally geodesic hypersurfaces 
corresponding to the complementary regions of the train track.
(Here the outer-normal is in the concave directions.)

%

We choose one of the normal vectors and a point $x_i$ on it 
$\gamma$-close to $I_i$, where $\gamma$ is a small positive number.
We push all the points of $I_i$ to $x_i$ 
to obtain a train track $\tau_{\eps, \delta, \gamma}$ 
and the complementary regions 
move accordingly to geodesic triangles with edges in 
the train track $\tau_{\eps, \delta, \gamma}$.


We claim that then the triangles are very close to the original 
triangles in their normal directions as well as in the
Hausdorff distance since the edge lengths of the triangles 
are bounded below.
Since there are no rapid turning 
of the complementary geodesic triangles, 
we can be assured that the new surface is imbedded
by integrations. 
This completes the proof of (i). 

The vertices corresponding to $I_1, \dots, I_n$ 
in the pleated open surfaces are s-vertices: 
The leaves of the laminations are moved in 
the normal direction which is a concave direction
Thus the leaf is bent toward the concave direction.
The triangles are of almost the same direction as before. 
Lemma \ref{lem:cross} implies the result. 
This completes the proof of (iii) and (iv).


(i) This matters for crescents that are inner if 
the perturbations are inner and ones that are outer if 
the perturbations are outer: In other cases, 
Lemma \ref{lem:crescentiso} shows that reversing 
the perturbation process gives us back all of our old
crescents preserving the $I$-part hypersurface.

Suppose that we have the perturbed sequence 
$\Sigma'_i$ closer and closer to $\Sigma$, 
and there exists a sequence of crescent $\calR_i$ 
for $\Sigma'_i$ not contained in
a certain neighborhood of the union of 
crescents for $\Sigma$. Then the limit $\calR$ of $\calR_i$
is still a crescent for $\Sigma$ and is not in 
the neighborhood. This is absurd.

\end{proof}


The move above described in Theorem \ref{thm:convexperturb} 
is said to be the {\em convex perturbation}.

\begin{cor}\label{cor:convperturb}
\begin{itemize} 
\item After the convex perturbation, the outer level of the resulting 
surface $\Sigma'$ is less or equal to the level of 
$\Sigma$ if we moved carefully, i.e., we do the small truncation moves 
first and then perturb vertices by sufficiently small amounts.
\item The inner level stays the same or become smaller than 
that of $\Sigma$. 
\item In particular, if $\Sigma$ has 
no crescent, then $\Sigma'$ contains no crescent
and hence $\Sigma'$ is s-imbedded triangulated surface.
Moreover, s-vertices that were at the boundary 
of the pleated part become strict s-vertices. 
\end{itemize}
\end{cor}
\begin{proof}
Let $n$ be the level of $\Sigma$. 
After the small truncation moves at (*) and (**), 
we obtain a surface $\tilde \Sigma''$ with level 
less than or equal to $n$. Also, the nature of 
truncation shows us that the union of 
crescents is the subset of that of $\tilde \Sigma$.
These all follow as in 
the proof of Proposition \ref{prop:perturb}.

In the pleated part, we define the {\em minimal pleated part} 
as the convex hull of the union of bi-infinite pleating 
geodesics with respect to the intrinsic metric  
obtained by piecing the pleated parts together. 
The minimal pleated part is a subsurface of 
the pleated part which is the closure of the union of 
totally geodesic subsurfaces bounded by 
the bi-infinite pleating geodesics. 

The minimal pleated part may meet the crescents 
only tangentially at the $I$-part. If not, then 
the minimal pleated part must pass through the $I$-part 
and this implies that a bi-infinite pleating 
geodesic pass through the $I$-part by the above 
paragraph. By Corollary \ref{cor:avoidg} (iii), 
this is a contradiction. 

We now move vertices of the train tracks of 
the pleating laminations by a very small amount 
according to Theorem \ref{thm:convexperturb}. 
The crescents do move in its $\eps$-neighborhood. 

Suppose now $n \geq 0$.
We choose the $\eps$-neighborhood sufficiently small 
so that any new component of $\tilde \Sigma$ intersected with 
crescents may not arise as $\tilde \Sigma$ is deformed: 
First, the small truncated places may be avoided 
by taking $\eps$-sufficiently small. 
Second, since any crescent during the perturbation cannot
meet the minimal pleated part in its interior, 
the $I$-part of crescents and crescents themselves
close to the minimal pleated part
lie below the minimal pleated part or may meet the minimal
pleated part but cannot pass through it
by Corollary \ref{cor:avoidg}.
Since the perturbation is in the concave direction,
our surface moved away from the crescents,
and outer crescents cannot achieve level $n+1$ 
or higher level.

Suppose that the inner level is $\geq 0$. 
Since our move is outward
we can push the $I$-parts of the inner crescents inward toward 
themselves and these are all the inner crescents obtainable
as before by Proposition \ref{prop:perturb}.
Therefore, the level can only decrease or stay the same
as after the crescent isotopy, i.e., before the convex
perturbation.

Suppose that the inner level was $-1$, i.e., there were no inner 
crescents. If $\Sigma'$ has an inner crescent, then 
since $\Sigma$ is obtainable by inner direction move from 
$\Sigma'$ by reversing our isotopy, we see that 
the inner crescents for $\Sigma$ exists by preserving
the I-parts. This is a contradiction. Thus, the inner 
level of $\Sigma'$ is $-1$. 

Suppose now $n = -1$. Then $\Sigma$ has no concave 
vertex. Thus the vertices at the boundary of the pleated 
part are all s-vertices. We move these
so that they become strict s-vertices: 
The boundary vertex lies on a simple closed curve in 
the boundary of the pleated part. Since there are no 
convex vertices, the curve is a concavely curved one with respect to 
the pleated part. At each vertex there exists a straight line 
through it in the pleated part.
The holonomy of the curve preserves a plane. We can push these vertices up 
in a parallel manner. Then the planes adjacent to 
the edges of the curve now begin to meet in angles $< \pi$. 
Thus, the straight lines get bent and the angles.
By a small perturbations, all moved vertices become 
strict s-vertices by Lemma \ref{lem:deforming}.
Therefore, all vertices of 
$\Sigma'$ are strictly s-vertices or convex vertices. 
Thus the outer level of $\Sigma'$ is $-1$.

\begin{figure}[ht]
\centerline{\epsfxsize=3.5in \epsfbox{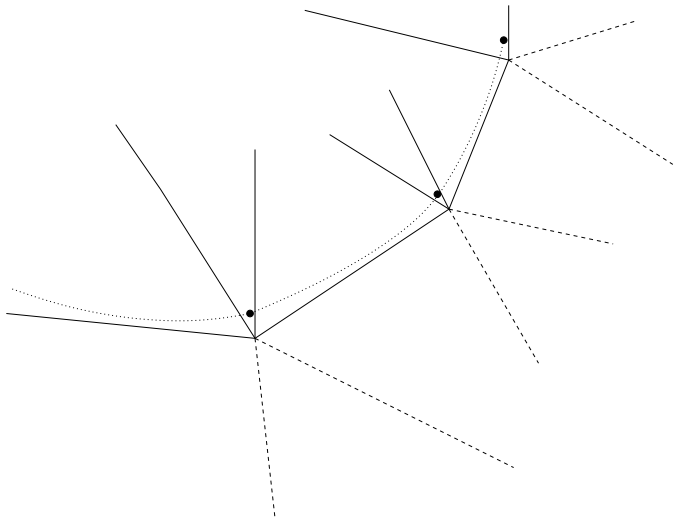}}
\caption{\label{fig:svertices} Making the boundary 
vertices of the pleated part into s-vertices.}
\end{figure}
\typeout{<>}

This implies that there are no outer crescents
since any crescents has an affine function $0$ on
the $I$-level and hence some extreme points under this function
in the $\alpha$-parts which gives us a concave vertex. 

If the inner level of $\Sigma$ is $-1$ also, 
then we see that $\Sigma'$ has no outer or inner crescents. 
\end{proof}


\section{Isotopy sequences}

In this section, we describe the process
of our modification of $\Sigma$ of highest-level $n$:
\begin{itemize}
\item[(i)] the small truncation moves for outer highest-level crescents.
\item[(ii)] the crescent 
isotopy for outer highest-level crescents. 
\item[(iii)] the convex perturbations which involves 
small truncation moves.
\end{itemize}
We do the same steps for the inner highest-level crescents, 
i.e., at the level $n$.
We do this until there are no crescents, i.e., to level $0$ 
inner and outer crescent moves and convex perturbations. 
The result is an isotopied $\tilde \Sigma$ with only 
s-vertices as there are no crescents
since a crescent must have a convex or concave vertex by
taking an extreme point for an affine function constant 
on the $I$-part.
Finally, we prove our main theorem. 


Also, the manifold bounded by the perturbed $\Sigma'$ still contains
$\cl_1, \dots, \cl_n$ nearby.



\begin{thm}\label{thm:aftermoves} 
Let $\Sigma$ have the highest-level equal to $n$. 
Let $\Sigma'$ be obtained from 
$\Sigma$ by 
\begin{itemize}
\item[(1)] small truncation moves for the level-$n$ outer crescents, 
\item[(2)] crescent-isotopy for level-$n$ outer crescents, and 
\item[(3)] convex perturbations. 
\end{itemize}
Then $\Sigma'$ is a triangulated surface:
\begin{itemize}
\item[(i)] the outer level of $\tilde \Sigma'$ 
is strictly less than that of $\tilde \Sigma$. 
\item[(ii)] In particular, if the outer level of $\tilde \Sigma$
was $0$, then $\tilde \Sigma'$ has no outer crescents. 
\item[(iii)] The inner level of $\tilde \Sigma'$ is 
less than or equal to that of $\tilde \Sigma$.
\end{itemize}
The statements also hold if 
every word ``outer'' were replaced by ``inner''
and vice versa. 
\end{thm} 
\begin{proof}  
(i) This was Corollary \ref{cor:convperturb}.

(ii) follows from (i). 

(iii) Given an inner crescent $\calR$ for the perturbed $\tilde \Sigma'$, 
we can isotopy its $\alpha$-parts to obtain an inner crescent 
for $\tilde \Sigma$ by Lemma \ref{lem:crescentiso} 
since the reverse movement is inward while so is $\calR$.

For $\calR$ to decrease its folding number, 
while $\tilde \Sigma'$ is isotopied back to $\tilde \Sigma$,
we need a convex vertex leaving $\calR$ through $I_{\calR}$. However, all of 
the above process reversed does not do this kind of movement.

The final statement is obvious.
\end{proof}

Let $n$ be the highest level for $\tilde \Sigma$.
Once, we do the outer highest-level crescent moves of level $n$, 
then we do the inner highest-level crescent moves of level $n$. 
By above Theorem \ref{thm:aftermoves}, 
we see that the inner and outer level strictly decreases
until there are no more crescents. 

Moreover the union of the set of crescents 
is contained in the $\eps$-neighborhood of 
the union of the set of crescents in the previous step. 

\begin{prop}\label{prop:geodesicavoid}
During the moves, the closed geodesic 
in the inner {\rm (}resp. outer {\rm )} component 
of complement of $\tilde \Sigma$ remains inside the $\eps$-neighborhood 
of the closure of the inner {\rm (}resp. outer {\rm)} 
component of $\tilde \Sigma'$
for arbitrarily small $\eps > 0$. 
\end{prop}
\begin{proof}
We use Corollary \ref{cor:avoidg}.
Only problem occurs when the $\partial_I$-part contains 
a lift $l$ of $\cl_i$. This implies $l \subset I_\calR$ 
for a secondary highest-level crescent $\calR$.
Hence, $l$ is homotopic into $\Sigma$.
Moreover, the holonomy $g$ along the closed geodesic has 
no rotation component since otherwise 
\[\bigcup_{i \in \bZ} g^i(\calR)\]
covers around $l$ to produce an annulus which 
compactifies to a sphere in $\hyps\cup \SI^2_\infty$.
This contradicts the assumption that $\Sigma$ is 
not homeomorphic to a torus.

After the crescent isotopy, $l$ is contained in 
a pleated-triangulated surface $\tilde \Sigma'$. 
We choose the convex perturbation that pushes $\tilde \Sigma'$ off $l$ 
by a small amount less than $\eps$. We find an arc $\alpha$ 
in pushed off $\tilde \Sigma''$ near $l$ with the same endpoints as $l$
and a neighborhood of $\alpha$ in $\tilde \Sigma''$ has s-vertices only.
Since the holonomy has no rotation component, we may assume 
that $\alpha$ is a convex arc nearly on a totally geodesic hypersurface $P$ 
passing through $l$. 
Any crescent during the isotopies and moves
that come later, it will either contain $l$ in its $I$-part or 
be a distance away. Draw a thin strip $A$ with boundary $l \cup \alpha$.
Then any innermost highest-level crescent $\calS$ meeting 
$\alpha$ and to-wards the direction of $l$ will meet $A$ but not 
contain any part of $l$ in its interior. Hence 
$\Lambda(\calS)$-isotopy
will not move $\tilde \Sigma''$ away from $l$.  

If $\calS$ is toward the opposite direction to $l$, 
$\alpha_\calS$ cannot meet $\alpha$ due to the local 
nature of $\alpha$: The totally geodesic hypersurface $P'$
containing $I_{\calS}$ meet $\alpha$ union its endpoints 
at least two points while the part of $\alpha$ 
in between must be in $\alpha_{\calS}$.
We see that $P \cap P'$ contains a geodesic $l'$ very close to $l$. 
This means $I_{\calS}$ contains $l'$. Since 
$l'$ is in the side of $l$, this is a contradiction. 

By induction, we can have that the end result surface 
$\tilde \Sigma$ is as close to $l$ as we want.
\end{proof}

Furthermore, the final perturbation gives 
us a triangulated surface isotopic to $\Sigma$. 


During the isotopy, we see that the support of the isotopy 
is included in the union of crescents. 

Applying the process in the beginning of the section, 
we can isotopy the boundary of $M$ to 
an s-surface. This follows since if there were 
a convex or concave vertex of $\partial M$ 
after all the moves, we find a crescent 
of level $0$. 
This completes the proof of Theorem C.

Here, $M'$ is isotopic to $M$ of course. 

\begin{rem}\label{rem:svertices} 
The vertices of the boundary of $M'$ are 
strict s-vertices: We start from general position 
$\partial M$ so that each s-vertex is a strict s-vertex
(see Theorem \ref{thm:convexperturb}).
Newly created vertices in the interior of the pleated 
part after the crescent isotopies are all strict s-vertices.
The boundary vertices of the pleated part at level $0$ 
are perturbed to be strict s-vertices
(see Corollary \ref{cor:convperturb}).
Finally, we no longer 
have any convex or concave vertex left at $M'$.
\end{rem}

\begin{rem}\label{rem:lamination} 
In fact, we can even assume that $\cl_i$ are 
geodesic laminations since we can modify 
Proposition \ref{prop:geodesicavoid}
using Proposition \ref{prop:avoidg}.
\end{rem}


\part{General hyperbolic $3$-manifolds and convex hulls of 
their cores}

\section{Introduction to Part 2}




The purpose of this paper is to show the significance of $2$-convexity.


A hyperbolic triangle is a subset of a metric space isometric 
with a hyperbolic triangle in the hyperbolic plane $\hypp$.
If the ambient space is a $3$-dimensional metric space, 
then we require it to be totally geodesic as well 
and develop into a totally geodesic plane in the hyperbolic 
space $\hyps$. 

A general hyperbolic manifold $M$ is a metric space 
with a locally isometric immersion $\dev: \tilde M \ra \hyps$
from its universal cover $\tilde M$. 
$\dev$ has an associated homomorphism 
$h: \pi_1(M) \ra \Isom(\hyps)$ given 
by $\dev \circ \gamma = h(\gamma)\circ \dev$ 
for each deck transformation $\gamma \in \pi_1(M)$. 
We require that the boundary of $M$ 
triangulated by totally geodesic hyperbolic triangles. 

A $2$-convex general hyperbolic manifold is 
a general hyperbolic manifold
so that a isometric imbedding from $T^o \cup F_1 \cup F_2 \cup F_3$ 
where $T^o$ is a hyperbolic tetrahedron in the standard 
hyperbolic space $\hyps$  and $F_1, F_2, F_3$ its faces extend to 
an imbedding from $T$.

\begin{thm}\label{thm:Gromovh}
The universal cover of 
a $2$-convex general hyperbolic manifold is Gromov hyperbolic. 
\end{thm}

An {\em h-map} from a triangulated hyperbolic surface is 
a map sending hyperbolic triangles to hyperbolic triangles 
and the sum of angles of image triangles around a vertex is 
greater than or equal to $2\pi$.

\begin{thm}\label{thm:GaussB2}
Let $\Sigma$ be a compact h-mapped surface relative to 
$v_1, \dots, v_n$ into a general hyperbolic manifold
and suppose that each arc in $\partial \Sigma - \{v_1, \dots, v_n\}$ 
is geodesic.
Let $\theta_i$ be the exterior angle of $v_i$ 
with respect to geodesics in the boundary of $\Sigma$.
Then 
\begin{equation}\label{eqn:GB2}
\Area(\Sigma) \leq \sum_i \theta_i - 2\pi \chi(\Sigma). 
\end{equation} 
\end{thm}

The convex hull of a homotopy-equivalent 
closed subset of a general hyperbolic manifold $M$
is the image in $M$ of the smallest closed convex subset of
containing the inverse image of the closed subset 
in the universal cover $\tilde M$ of $M$.
We can always choose the core 
$\core$ to be a subset of $M^o$.

\begin{thm}\label{thm:convh}
Let $\convh(\core)$ be the convex hull of the core $\core$ of 
a $2$-convex general hyperbolic manifold.
We assume that $\core$ is chosen to be a subset of 
$M^o$ and $\partial \core$ is s-imbedded. 
Suppose that $\convh(\core)$ is compact.
Then $\convh(\core)$ is homotopy equivalent to the core and 
the boundary is a truly pleated-triangulated h-surface.
\end{thm}


In the preliminary, Section 8, 
we recall the definition of 
\cat$\kappa$-spaces for $\kappa \in \bR$
using geodesics and triangles. 
We also define $M_\kappa$-spaces, the simplicial metric 
spaces needed in this paper. We discuss the link conditions
to check when $M_\kappa$-space is \cat($\kappa)$-space, 
the Cartan-Hadamard theorem, and Gromov boundaries of these spaces.
Next, we discuss the $2$-dimensional versions
of these spaces. Define the interior angles, Gauss-Bonnet 
theorem. Finally, we discuss general hyperbolic manifolds. 

In Section 9, we show that the universal cover of 
a $2$-convex general hyperbolic 
manifold, which we used a lot in Part 1, is 
a $M_{-1}$-simplicial metric space and a \cat($-1$)-space
and a visibility manifold. 
Next, we define h-maps of surfaces. These are similar 
to hyperbolic surfaces as defined by Bonahon, Canary, and Minsky. 
We define A-nets, a generalization of triangles,
and $\alpha$-nets a simplicial approximation and a h-map to
A-nets. We show that maps from surfaces can be homotopied to
h-maps relative to a collection of boundary points
in $2$-convex general hyperbolic manifolds.
We prove the Gauss-Bonnet theorems for such surfaces 
and find area bounds for polygons.

In Section 10, we discuss the convex hull of 
the core $\core$ in a general hyperbolic manifold.
First, we show that the convex hull and $\core$ is homotopy 
equivalent. We finally show that 
the boundary of the convex hull can be deformed to 
a nearby h-imbedded 
surface, which is truly pleated-triangulated.
We show this by finding geodesic in the boundary 
through each point of the boundary.



\section{Hyperbolic metric spaces}

\subsection{Metric spaces, geodesics spaces, and cat(-1)-spaces.}
\label{sec:metricspaces}

We will follows Bridson-Haefliger \cite{BH}.

Let $(X, d)$ be a metric space. A {\em geodesic path} from 
a point $x$ to $y$, $x, y \in X$ is a map 
$c: [0, l] \ra X$ such that $c(0)=x$ and $c(l) = y$ 
and $d(c(t), c(t')) = |t - t'|$ for all $t, t' \in [0, l]$. 

A {\em local geodesic} is a map $c: I \ra X$ from an interval 
$I$ with the property that for every $t \in I$ there exists $\eps > 0$ 
such that $d(c(t'), c(t'')) = |t'- t''|$ for 
$t', t''$ in the $\eps$-neighborhood of $t$ in $I$. 

$(X,d)$ is a {\em geodesic metric space} if every pair of points of 
$X$ is joined by a geodesic. 

We denote by $\bE^2$ the plan $\bR^2$ with the standard Euclidean metric. 
A {\em comparison triangle} in $\bE^2$ 
of a triple of points $(p, q, r)$ in $X$
is a triangle in $\bE^2$ with vertices $\bar p, \bar q, \bar r$ 
such that $d(p,q) = d(\bar p, \bar q)$, 
$d(q, r) = d(\bar q, \bar r)$, and 
$d(r, p) = d(\bar r, \bar p)$. 
This is unique up to isometries of $\bE^2$. 

The interior angle of the comparison triangle at $\bar p$ 
is called the comparison angle between $q$ and $r$ at $p$ 
and is denoted $\bar \angle_p(q, r)$. 

Let $c:[0, a] \ra X$ and $c': [0, b] \ra X$ be two 
geodesics with $c(0) = c'(0)$. 
We define the {\em upper angle} $\angle_{c, c'} \in [0, \pi]$ 
between $c$ and $c'$ to be 
\begin{equation}\label{eqn:angle} 
\angle(c, c') := \limsup_{t, t' \ra 0} \bar \angle_{c(0)}(c(t), c'(t'))
\end{equation}
The angle exists in strict sense if the limsup equals the limit. 

The angles are always less than or equal to $\pi$ by our construction. 
We define angles greater than $\pi$ in 
two-dimensional spaces by specifying sides 
and dividing the side into many parts
(see Subsection \ref{subsec:singhyp}).

We say that a sequence of closed subsets $\{K_i\}$ converge to 
a closed subset $K$ if for any compact subset $A$ of $X$, 
$\{A \cap K_i\}$ converges to $A \cap K$ in Hausdorff sense.

Let $(X, d)$ be a metric. We can define a {\em length-metric} $\bar d$ so that 
$\bar d(x, y)$ for $x, y \in X$ is defined as the infimum of the lengths of 
all rectifiable curves joining $x$ and $y$. 
We note that $d \leq \bar d$ and 
$(X, d)$ is said to be a length space if $\bar d = d$. 

\begin{prop}[Hopf-Rinow Theorem] 
Let $(X, d)$ be a length space. If $X$ is complete and locally 
compact, then every closed bounded subset of $X$ is compact 
and $X$ is a geodesic space. 
\end{prop} 

As an example, a Riemannian space with path-metric is a geodesic 
metric space. A covering space of a length space has an obvious induced 
length metric



We define $M_\kappa$ to be the $3$-sphere of constant curvature $k$, 
Euclidean space, or the hyperbolic $3$-space of constant curvature $k$
depending on whether $\kappa > 0,=0,< 0$ respectively. 

Let $D_\kappa$ denote the diameter of $M_\kappa$ if $\kappa > 0$ 
and let $D_\kappa = \infty$ otherwise. 

Let $(X, d)$ be a metric space. Let 
$\tri$ be a geodesic triangle in $X$ with parameter 
less than $2D_\kappa$  and $\bar \tri$ the comparison
triangle in $M_\kappa$. 
Then $\tri$ is said to satisfy {\em \cat($\kappa$)-inequality} if 
$d(x, y) \leq d(\bar x, \bar y)$ for 
all $x, y$ in the edges of $\tri$ and their comparison 
points $\bar x, \bar y$, i.e., of same distance from 
the vertices, in $\bar \tri$.
If $\kappa < 0$, a {\em \cat($\kappa$)-space} is a geodesic space all of whose 
triangles bounded by geodesics satisfy \cat($\kappa$)-inequality. 
If $\kappa > 0$, then $X$ is called a {\em \cat($\kappa$)-space} if 
$X$ is $D_\kappa$ geodesic and all geodesic triangles in 
$X$ of perimeter less than $2D_\kappa$ satisfy 
the \cat($\kappa$)-inequality. 

Angles exists in the strict sense for 
\cat($\kappa$)-spaces if $\kappa \leq 0$. 

A \cat($\kappa$)-space is a \cat($\kappa'$)-space if $\kappa \leq \kappa'$.

A \cat($0$)-space $X$ has a metric $d:X \times X \ra \bR$
that is convex; i.e., given any two geodesics $c:[0,1] \ra X$ 
and $c':[0,1] \ra X$ parameterized proportional to length, 
we have 
\begin{equation}\label{eqn:geoconv} 
d(c(t), c'(t)) \leq (1-t) d(c(0), c'(0)) + td(c(1), c'(1)).
\end{equation} 


A geodesic $n$-simplex in $M_\kappa$ is the convex hull of 
$n+1$ points in general position. 

An {\em $M_\kappa$-simplicial complex} $K$ is 
defined to be the quotient space of the 
disjoint union $X$ of a family of geodesic $n$-simplicies 
so that the projection $q: X \ra K$ induces the injective projection
$p_\lambda$ for each simplex $\lambda$
and if $p_\lambda(\lambda) \cap p_{\lambda'}(\lambda')\ne \emp$, 
there exists an isometry $h_{\lambda, \lambda'}$ from a face of 
$\lambda$ to $\lambda'$ such that 
$p_\lambda(x) = p_{\lambda'}(x')$ if and only if 
$x' = h_{\lambda, \lambda'}(x)$. 

In this paper, we will restrict to the case when locally there 
are only finitely many simplicies, i.e., $X$ is locally convex. 
We do not assume that we have a finite isometry types of 
simplicies as Bridson does in \cite{Bridson}.

A geodesic link of $x$ in $K$, denoted by $L(x, K)$  is the set of 
directions into the union of simplicies containing $x$. 
The metric on it is defined in terms of angles. 
(For details, see Chapter I.7 of \cite{BH}.)

\begin{defn}\label{defn:linkcond} 
An $M_\kappa$-simplicial complex satisfies the link condition if 
for every vertex $v$ in $K$, the link complex 
$L(v, K)$ is a \cat($1$)-space. 
\end{defn}

The following theorem can be found in Bridson-Haefliger \cite{BH}:
\begin{thm}[Ballman]\label{thm:linkcond} 
Let $K$ be a locally compact $M_\kappa$-simplicial complex.
If $\kappa \leq 0$, then the following conditions are equivalent\rmc
\begin{itemize} 
\item[(i)] $K$ is a \cat($\kappa$)-space. 
\item[(ii)] $K$ is uniquely geodesic. 
\item[(iii)] $K$ satisfies the link condition and contains no isometrically
imbedded circle. 
\item[(iv)] $K$ is simply connected and satisfies the link condition. 
\end{itemize} 
If $\kappa > 0$, then the following conditions are equivalent \rmc
\begin{itemize}
\item[(v)] $K$ is a \cat($\kappa$)-space. 
\item[(vi)] $K$ is $\pi/\sqrt{\kappa}$-uniquely geodesic. 
\item[(vii)] $K$ satisfies the link condition and contains no
isometrically embedded circles of length less than $2\pi/\sqrt{\kappa}$. 
\end{itemize} 
\end{thm} 
\begin{proof} 
See Ballmann \cite{Ballmann} or Bridson-Haefliger \cite{BH}.
\end{proof}

\begin{lem}\label{lem:2Dlinkcond} 
A $2$-dimensional $M_\kappa$-complex $K$ satisfies the link condition 
if and only if for each vertex $v \in K$, every injective loop 
in $Lk(v, K)$ has length at least $2\pi$. 
\end{lem} 


\begin{defn} 
A metric space $X$ is said to be of {\em curvature $\leq \kappa$} 
if it is locally isometric to a $\cat(\kappa)$-space, 
i.e., for each point $x$ of $X$, 
there exists a ball which is a $\cat(\kappa)$-space. 
\end{defn}

\begin{thm}[Cartan-Hadamard]\label{thm:CaHa}
Let $X$ be a complete metric space. 
\begin{itemize} 
\item[(i)] If the metric on $X$ is locally convex, then 
the induced length metric on the universal cover $\tilde X$
is globally convex. 
{\rm (}In particular, there is a unique geodesic connecting 
two points of $\tilde X$, and geodesic segments in $\tilde X$ vary 
continuously with respect to their endpoints.{\rm )} 
\item[(ii)] If $X$ is of curvature $\leq \kappa$ where $\kappa \leq 0$, 
then $\tilde X$ is a \cat$(\kappa)$-space. 
\end{itemize} 
\end{thm}

Let $\delta$ be a positive real number. 
A geodesic triangle in a metric space 
$X$ is said to be {\em $\delta$-slim} if each of its sides is 
contained in the $\delta$-neighborhood of the union of 
the other two sides. 

For $\kappa < 0$, \cat($\kappa$)-space is $\delta$-hyperbolic.

For positive real numbers $\lambda, \eps$, 
{\em $(\lambda, \eps)$-quasi-geodesic} in $X$ is a map 
$c: I \ra X$ such that 
\begin{equation}
1/\lambda |t - t'| - \eps 
\leq d(c(t), c(t')) \leq \lambda |t - t'| + \eps 
\end{equation} 

Let $X$ be a $\delta$-hyperbolic space. 
Two quasi-geodesic rays $c, c'$ are 
{\em equivalent} or {\em asymptotic} if 
their Hausdorff distance is finite, or, equivalently
$\sup_t d(c(t), c'(t))$ is finite.
We define the {\em Gromov boundary} $\partial X$  
as the space of equivalence classes of quasi-geodesic
rays in $X$. One can show that $\partial X$ is 
the space of equivalence classes of geodesics rays 
as well. 

If $X$ is a proper metric space, then $X$ is a visibility 
space: For each pair of points $x$ and $y$ in $\partial X$, 
there exists a geodesic limiting to $x$ and $y$. 
Topology and metrics are given on $\partial X$ to 
compactify $X \cup \partial X$. 
The group of isometry acts as homeomorphisms on $\partial X$.

\subsection{Sigular hyperbolic surfaces}\label{subsec:singhyp} 


A {\em hyperbolic triangle} in a metric space is 
a subset isometric to a triangle in $\hypp$ bounded 
by geodesics. 
Sometimes, we need a {\em degenerate hyperbolic triangle}. 
It is defined to be a straight geodesic segment or a point
where the vertices are defined to be the two endpoints and 
a point, which may coincide. 

A {\em hyperbolic tetrahedron} in a metric space is a subset 
isometric to a tetrahedron in $\hyps$ bounded by 
four totally geodesic planes with six edges geodesic 
segments and four vertices. 
Again degenerate ones can obviously be defined on a hyperbolic triangle, 
a segment, and a point with various vertex and edge structures.

A {\em hyperbolic cone-neighborhood} of a point $x$ in a surface 
$\Sigma$ with a metric is a neighborhood of $x$ which divides into 
hyperbolic triangles with vertices at $x$.
The {\em cone-angle} is the sum of angles of the triangles 
at $x$.  
The set of singular points is denoted by 
$\sing(\Sigma)$ and the cone-angle at $x$ by 
$\theta(x)$. 

By a {\em singular hyperbolic surface}, we mean 
a complete metric space $X$ locally isomorphic to 
a hyperbolic plane or a hyperbolic cone-neighborhood 
with cone-angle $\geq 2\pi$ so that the set of 
singular points are discrete.  
We will also require that $X$ is triangulated by 
hyperbolic triangles in this paper (i.e., is 
a metric simplicial complex 
in the terminology of \cite{BH}). 

By Lemma \ref{lem:2Dlinkcond}, the universal cover $\tilde X$ of $X$ 
is a \cat($-1$)-space. 

Let $X$ be a singular hyperbolic surface.
Clearly, $X$ has an induced length metric
and is a geodesic space. 

We say that a geodesic is {\em straight} if it is 
a continuation of geodesics in hyperbolic triangles 
meeting each other at $\pi$-angles in the intrinsic sense.

We can measure angles greater than $\pi$ in 
singular hyperbolic surface by dividing the 
angle into smaller ones. In this case, we need 
to specify which side you are working on. 

In general a path is {\em geodesic} if it is a continuation  
of straight geodesic meeting each other at 
greater than or equal to $\pi$-angles from both sides. 

We also say that a boundary point $x$ is {\em bent} if 
the two straight geodesics end at the point
not at $\pi$-angle in the interior. 
We define $\theta(x)$ to be $\pi$ minus the interior angle. 
It could be negative. 
We will denote by $\sing(\partial X)$ the 
set of bending points. 


\begin{prop}[Gauss-Bonnet Theorem] \label{prop:GaussBonnet}
Let $\Sigma$ be a compact singular hyperbolic surface
with piecewise straight geodesic boundary. 
Then 
\begin{equation}\label{eqn:GaussBonnet} 
-\Area(\Sigma) + \sum_{v \in \sing(\Sigma)} (2\pi-\theta(v)) 
+ \sum_{v \in \sing(\partial \Sigma)} \theta(v) 
= 2\pi \chi(\Sigma).
\end{equation}
\end{prop}

From the Gauss-Bonnet theorem, we can show that
there exists no disk bounded by two geodesics. 
This follows since if such a disk exists, 
then $\theta(v) \geq 2\pi$ for all singular, 
the exterior angles at virtual vertices $\leq 0$,
the exterior angles at common end points $< \pi$,
and the area is less than $0$. 

This implies: 
Given a compact singular hyperbolic space
and a closed curve, we can homotopy the curve into
a closed geodesic, and the closed geodesic is 
unique in its homotopy class. 


Moreover, two closed geodesics meet in a minimal number of 
times up to arbitrarily small perturbations: 
that is, the minimum of geometric intersection number 
under small perturbation is the true minimum under 
all perturbations. 
(We may have two geodesics agreeing on an interval 
and diverging afterward unlike the hyperbolic plane.)

\subsection{ General hyperbolic $3$-manifolds} 

By a {\em general hyperbolic manifold}, 
we mean a manifold $M$ with an atlas of charts to 
$\hyps$ with transition maps in $\Isom(\hyps)$. 
The metric on it will be the length metric given 
by the induced Riemannian metric. We require the metric 
to be complete. As a consequence, this is a geodesic space
by local compactness \cite{BH}. 
In general we assume that $\partial M$ is not empty. 
If it is not geodesically complete, 
$M$ need not be a quotient of $\hyps$ which 
are the usual subject of the study in $3$-manifold theory.

Also, we will require general hyperbolic manifolds to have 
hyperbolic triangulations, i.e., it has a triangulation 
so that each tetrahedron is isometric with 
a hyperbolic tetrahedron in $\hyps$. 
Moreover, we assume that the vertices of 
the triangulations are discrete and
the induced triangulation on the universal cover 
map under $\dev$ to a collection of tetrahedra in 
general position in $\hyps$. 
We also require the following mild condition: 
Every boundary point of a general hyperbolic 
manifold has a neighborhood isometric with a subspace of 
a metric-ball in $\hyps$. By subdivisions and small modifications, we can 
always achieve this condition. 

We will say that $M$ is {\em locally convex} if there is
an atlas of charts where chart images are convex subsets of 
$\hyps$. Thus, $M$ is locally convex if $\partial M$ is empty. 
(In this paper, we will be interested in the 
non-locally-convex manifolds.)

Given a general hyperbolic manifold $M$, 
its universal cover $\tilde M$ has an immersion
$\dev: \tilde M \ra \hyps$, which is not in general 
an imbedding or a covering map, and a homomorphism $h$ 
from the deck transformation group $\pi_1(M)$ to 
$\Isom(\hyps)$ satisfying 
\begin{equation}\label{equ:dev} 
\dev\circ \vth = h(\vth) \circ \dev, \vth \in \pi_1(M). 
\end{equation}
$\dev$ is said to be a {\em developing map} and 
$h$ a holonomy homomorphism. 

\begin{thm}[Thurston]\label{thm:localconv} 
Let $M$ be a metrically complete $3$-manifold and is 
locally convex. Then its developing map $\dev$ is 
an imbedding onto a convex domain, and $M$ is 
isometric with a quotient of a convex domain in $\hyps$
by an action of a Kleinian group.
\end{thm}
\begin{proof} 
See \cite{Thnote}.  
\end{proof} 

In this paper, we will often meet
{\em drilled hyperbolic manifolds} obtained 
by removing the interior of 
a codimenion-zero submanifold of a general hyperbolic manifolds. 
They are of course general hyperbolic manifolds.

Of course, special hyperbolic manifolds are general hyperbolic 
manifolds and drilled hyperbolic manifolds.



Since a general hyperbolic manifold has a geodesic metric, 
we can define geodesics. 
A {\em straight} geodesics is a geodesic which 
maps to geodesic in $\hyps$ under the charts. 
Geodesics are in general a union of straight geodesics. 
Thus, it has many bent points in general. 
The bent points in the interior of the geodesic segments 
are said to be {\em virtual vertices}.

We define angles as above for metric spaces. 
Then at a bent virtual vertex, the angle is
equal to $\pi$ since if not, then we can 
shorten the geodesics.

\begin{prop}\label{prop:2Dangle} 
Let $l$ be a geodesic with a bent point $x$ in 
its interior. Let $S$ be a simplicially immersed 
surface containing $l$
in its boundary. Give $S$ an induced length metric. 
Then the interior angle at $x$ in $S$ is always greater than 
or equal to $\pi$. 
\end{prop}
\begin{proof} 
If the angle is less than $\pi$, 
we can shorten the geodesic. 
\end{proof} 

Given an oriented geodesic $l_1$ ending at a point $x$ 
and an oriented geodesic $l_2$ starting from $x$, 
we can define an {\em exterior angle} between $l_1$ and $l_2$ to be 
$\pi$ minus the angle between the geodesic $l'_1$ 
with reversed orientation and the other one $l_2$. 



\section{ $2$-convex general hyperbolic manifolds 
and {h}-maps of surfaces}

\subsection{$2$-convexity and general hyperbolic manifolds}


In Part 1, we showed that a general hyperbolic 
manifold was $2$-convex if the vertices of the boundary 
were either s-vertices or convex vertices.


We recall the definition of $2$-convexity: 
\begin{defn}\label{defn:2conv}
A general hyperbolic manifold is 
{\em $2$-convex} if given a compact subset $K$ mapping to 
a union of three sides of a tetrahedron $T$ in $\hyps$ 
under a chart $\phi$ of atlas, there exists a subset $T'$ 
mapping to $T$ by a chart extending $\phi$.
\end{defn} 

\begin{prop}\label{prop:uni2conv} 
If $M$ is a $2$-convex general hyperbolic manifold, 
then $M$ is a $K(\pi_1(M))$, i.e., 
its universal cover is contractible.  
\end{prop} 
\begin{proof}
Since the universal cover $\tilde M$ has an affine 
structure with trivial holonomy induced from 
the affine space containing $\hyps$ from the Klein model, 
this follows from \cite{uaf}. Also,
this follows from Theorem \ref{thm:2convcat}. 
\end{proof}


\begin{thm}\label{thm:2convcat} 
Let $M$ be a $2$-convex general hyperbolic manifold. 
Then its universal cover $\tilde M$ is 
a $M_{-1}$-simplicial complex and a \cat($-1$)-space. 
($M$ has a curvature $\leq -1$.)
\end{thm} 
\begin{proof} 
Using Theorem \ref{thm:linkcond} (iv), we need to show that 
for each vertex $x$ of $\tilde M$, the link $P=L(x, \tilde M)$ is 
a \cat($1$)-space. 

To show $P$ is a \cat($1$)-space, we use (vii) of Theorem \ref{thm:linkcond};
i.e., we show that $P$ satisfies the link condition and contains 
no isometric circle of length $< 2\pi$. By the boundary condition on $M$, 
$P$ is isometric to the unit sphere if $x$ is the interior point
or is isometric to a subspace of the unit sphere if $x$ is 
the boundary point. Clearly the former satisfy 2-dimensional link conditions. 

Let $P$ be a proper subspace of a unit sphere $\SI^2$ and $c$ an isometrically 
imbedded circle of length $< 2\pi$. By Lemma \ref{lem:lengthhemisphere}, 
$c$ is disjoint from a closed hemisphere $H$ in $\SI^2$. 

The closed curve $c$ meets $\partial P$ since otherwise $c$ has 
to be a great circle of length $2\pi$ being a geodesic. 
However, $c$ may never cross the circle
$\partial P$ over. Let $D^1_c$ and $D^2_c$ denote the disks 
in $\SI^2$ bounded by $c$. Then $\partial P$ is a subset of 
$D^1_c$ or $D^2_c$. Assume without loss of generality that 
the former is true. Suppose that $H$ is a subset of $D^2_c$. 
Then $H \subset P$. Looking this from $x$ again, we see that 
$2$-convexity is violated. 

Suppose that $H$ is a subset of $D^1_c$. Let $H'$ be the complementary 
open hemisphere. Then $c \subset H'$ and $\delta P$ is 
a outside the disk $D^2_c$ in $H'$ bounded by $c$. 
Since $H'$ has a natural affine structure, let $K$ be 
the convex hull of $c$ in $K$. Let $y$ be an extreme point of 
$K$. Then $y \in c$ as well. Near $y$ inside $c$ there are no 
points of $\delta P$. Thus, we can shorten $c$ 
contrary to the fact that $c$ is a geodesic. 
\end{proof} 

\begin{lem}\label{lem:lengthhemisphere} 
Let $\gamma$ be a broken geodesic loop in the sphere $\SI^2$ of 
radius $1$. If the length of $\gamma$ is less than $2\pi$, 
then there exists an open hemisphere containing it 
{\rm (}and hence a disjoint closed hemisphere{\rm ).}
\end{lem}
\begin{proof} 
We can shorten the loop without increasing the number of 
broken points to a loop as short as we want. 
A sufficiently short loop is contained in an open hemisphere. 

Let $l_t, t\in [0,1]$ be a homotopy so that 
$l_1$ is the original loop and $l_0$ is a constant loop. 
Then let $A$ be the maximal connected set containing $0$ 
so that $l_t$ for $t \in A$ is contained in an open hemisphere, say $H_t$. 

The set $A$ is an open set since the small change in $l_t$ does not 
violate the condition. Suppose that the complement of $A$ is not empty. 
Let $t_0$ be the greatest lower bound of 
the complement of $A$. 
Then $l_{t_0}$ is contained in a closed hemisphere, say $H'$, since 
we can find a geometric limit of the closure of $H_t$s. 

Suppose that $\partial H' \cap l_{t_0}$ is contained in 
a subset of length strictly less than $\pi$. 
Then we can rotate $H'$ along a pivoting antipodal 
pair of points on $\partial H'$ outside the subset. 
Then the new hemisphere contains $l_{t_0}$ in its interior, 
a contradiction. 

Suppose that $\partial H' \cap l_{t_0}$ contains an antipodal points. 
Let $s_1$ and $s_2$ be the corresponding points of $[0,1]$
and suppose $0 < s_1< s_2 < 1$ without loss of generality. 
Then two arcs $l_{t_0}|[s_1, s_2]$ and $l_{t_0}| [s_2, 1] \cup [0, s_1]$,
must have length greater than or equal to $\pi$, a contradiction. 

Therefore, no subsegment in $\delta H'$ of length $\leq \pi$ contains 
$\partial H' \cap l_{t_0}$, and we may assume without loss of 
generality that there are three points $p_1, p_2, p_3$
in $\delta H' \cap l_{t_0}$ are not contained in 
a subsegment in $\delta H'$ of length $\leq \pi$ 
and no pair of them are antipodal. 

The sum of lengths of segments $\ovl{p_1p_2}, \ovl{p_2p_3}, \ovl{p_3p_1}$ 
equals $2\pi$. 
This is clearly less than or equal to that of $l_{t_0}$ since 
the shortest arcs are these segments. 
This is again a contradiction. 

Thus $A$ must be all of $[0,1]$. 
\end{proof}

The following proves Theorem \ref{thm:Gromovh} in detail.
\begin{prop}\label{prop:Mgeodcond} 
Let $\tilde M$ be a universal cover of 
a compact $2$-convex general hyperbolic 
manifold $M$. Then the following hold{\rm :} 
\begin{itemize}
\item $\tilde M$ is uniquely geodesic.
\item Geodesic segments of $\tilde M$ depend continuously on their endpoints. 
\item The metric is locally convex. 
\item $\tilde M$ is $\delta$-hyperbolic and hence it is a visibility manifold. 
\item $M$ has curvature $\leq -1$. 
\item Given any path class on $M$, there exists a unique geodesic 
segment, which depends continuously on endpoints.
\end{itemize}
\end{prop}
\begin{proof} 
These are direct consequences of the fact that 
$\tilde M$ is a \cat$(-1)$-space. 
\end{proof}

\subsection{h-maps of surfaces into $2$-convex general
hyperbolic manifolds}

A {\em triangulated hyperbolic surface} is a metric surface with 
or without boundary triangulated and each triangle is 
isometric with a hyperbolic triangle or a degenerate 
hyperbolic triangle in $\hypp$. 
A {\em half-space} of $\hyps$ is a subset bounded by 
a totally geodesic plane.  

\begin{defn}\label{defn:h-map} 
Let $\Sigma$ be a compact triangulated hyperbolic surface, $M$ 
a general hyperbolic $3$-manifold, and $f:\Sigma \ra M$
a map which sends each triangle to a hyperbolic triangle in $M$. 
Let $\partial \Sigma$ have distinguished vertices 
$v_1, \dots, v_n$. 
Then $f$ is an {\em h-map relative to} $\{v_1, \dots, v_n\}$ 
if the sum of the angles of the image triangles of 
the stellar neighborhood of each interior vertex 
$v$ is greater than or equal to $2\pi$ and 
the sum of angles of the image triangles of 
the stellar neighborhood of the boundary vertex $v$, $v \ne v_i$, 
is greater than or equal to $\pi$. 
\end{defn}  

An h-map is a completely analogous concept to a hyperbolic-map by
Bonahon \cite{Bonahon}, Canary and Minsky and so on.
Note that if the boundary portion between 
$v_i$ and $v_{i+1}$ is geodesic for each $i$, 
then the boundary angle conditions are satisfied also.



\begin{defn}\label{defn:A-net}
Given an arc or a point $\alpha$ and an arc $\beta$ in $M$, 
an {\em Alexandrov net} or {\em A-net with ends $\alpha$ and $\beta$}
is a map $f: I \times I \ra M$
so that $s \in I \mapsto f(t, s)$ is geodesic for 
each $s$ and $t \mapsto f(t,0)$ is $\alpha$ and 
$t \mapsto f(t, 1) \in \beta$.
\end{defn}

\begin{lem}\label{lem:closegeo} 
Let $M$ be a $2$-convex general hyperbolic manifold.  
Let $\gamma$ be a geodesic in $M$. Then for any geodesic $\gamma'$
sufficiently close to $\gamma$,
there exists a homotopy $H: I \times I \ra M$
so that the following hold: 
\begin{itemize}
\item $s \mapsto H(0,s)$ is $\gamma$ 
and $s \mapsto H(1, s)$ is $\gamma'$. 
\item $H$ is a simplicial map with 
a triangulation of $I\times I$ with 
all vertices at $\{0,1\} \times I$. 
\end{itemize}
\end{lem} 
\begin{proof} 
For each virtual vertex of $\gamma$, we choose 
$\eps> 0$ such that the $\eps$-neighborhood of 
the vertex is a stellar neighborhood. 
If $\gamma'$ is in an $\eps$-neighborhood of $\gamma$ 
for $\eps>0$ for any such $\eps$s, 
then we can find the desired $H$. 
\end{proof}

\begin{def}\label{defn:a-net} 
Given a point or an arc $\alpha$ and another arc $\beta$,
an {\em $\alpha$-net} $f:I \times I \ra M$ with ends $\alpha$ and 
$\beta$ is a map such that 
\begin{itemize}
\item $s \mapsto f(t_i, s)$ for a finite subset 
$\{t_1=0, t_2, \dots, t_n=1\}$ of $I$ is a geodesic 
for each $i$,
\item $t \mapsto f(t, 0)$ is $\alpha$ 
and $t \mapsto f(t, 1)$ is $\beta$. 
\item $f$ is an h-map relative 
to vertices of the arcs $\alpha$ and $\beta$
with a triangulation of $I\times I$ with all the vertices in 
$\{t_1, \dots, t_n\} \times I$.   
\end{itemize}
\end{def}

\begin{prop}\label{prop:a-netexist} 
Given a point or an arc $\alpha$ and another arc $\beta$,
there exists an $\alpha$-net with ends $\alpha$ and $\beta$.
\end{prop}
\begin{proof}
We find an A-net $f: I \times I \ra M$ 
with ends $\alpha$ and $\beta$. 
We take sufficiently many $t_i$'s so that 
geodesics $s \mapsto f(t_i, s)$ are very close. 
By Lemma \ref{lem:closegeo}, we can 
find a simplicial map $F: I \times I \ra M$. 
Since $s \mapsto F(t_i, s) = f(t_i, s)$ are 
geodesics, the sum of angles at each of 
the sides of a vertex on this geodesic 
is greater than $\pi$. Hence, the sum of 
angles at an interior vertex is greater than or 
equal to $2\pi$. At the vertices of $s \mapsto F(0, s)$ 
or $s \mapsto F(1, s)$, the sum of angles are 
greater than $\pi$. Therefore, $F$ is an h-map.
\end{proof} 


\begin{prop}\label{prop:himmexists} 
Let $\Sigma$ be a compact triangulated surface, $M$ a general 
hyperbolic $3$-manifold, and let $f:\Sigma \ra M$ be a map
with an injective induced homomorphism $f_*: \pi_1(\Sigma) \ra \pi_1(M)$. 
\begin{itemize}
\item Let $v_1, \dots, v_n$ be the distinguished vertices in $\partial \Sigma$
and $l$ be a union of disjoint simple closed curve in $\Sigma$
which is disjoint from $\{v_1, \dots, v_n\}$ 
and is a component of $\partial \Sigma$ or is disjoint from $\partial \Sigma$.
\item We suppose that $\{v_1,\dots, v_n\} \cup l$ is not empty.  
Suppose that $f$ maps each arc in $\partial \Sigma$ 
connecting two distinguished vertices to a geodesic
and each component of $l$ or $\partial \Sigma$ without 
any of $v_1, \dots, v_n$ to a closed geodesic. 
\end{itemize} 
Then in the relative homotopy class of $f$ with 
$f|\partial \Sigma$ fixed, 
there exists an h-map $f': \Sigma' \ra M$
relative to $v_1, \dots, v_n$ where $\Sigma'$ is 
$\Sigma$ with a different triangulation in general and 
$f'$ agrees with $f$ on $\partial \Sigma \cup l$. 
\end{prop}

From now on, we will just use $v_i$ for $f(v_i)$ and so on
since the reader can easily recognize the difference. 
By the angle of a triangle,
we mean the corresponding angle measured in the image triangle of $f$. 

\begin{proof} 
First, we find a topological triangulation $\Sigma$ 
so that all the vertices are in the union of 
$\{v_1, \dots, v_n\} \cup l \cup \partial \Sigma$. 
We find geodesics in the right path-class 
for each of the edges of the triangulations.
For each triangle, we extend by choosing 
a vertex and the opposite geodesic edges 
and finding $\alpha$-nets with these ends. 

At each interior point of an edge, 
we see that the sum of angles of any of its side is greater 
than or equal to $\pi$ since the edge is geodesic. 
Since $\alpha$-nets are h-maps, we see that 
the whole map is an h-map. 
\end{proof}

\subsection{Gauss-Bonnet theorem for h-maps}

\begin{prop}\label{prop:GB} 
Let $\Sigma$ be a compact h-mapped surface relative to $v_1, \dots, v_n$. 
Let $\theta_i$ be the exterior angle of $v_i$ 
with respect to geodesics in the boundary of $\Sigma$.
Then 
\begin{equation}\label{eqn:GB} 
\Area(\Sigma) \leq \sum_i \theta_i - 2\pi \chi(\Sigma). 
\end{equation} 
\end{prop}
\begin{proof} 
The interior angle with respect to $\Sigma$ is larger than 
the angle in $M$ itself. 
Thus the exterior angle with respect to $\Sigma$ is 
smaller than the exterior angle in $M$. 

Since the interior vertices have the angle sums greater than or 
equal to $2\pi$ and the boundary virtual vertices have 
the angle sum greater than or equal to $\pi$, 
the proposition follows. 
\end{proof}

An $n$-gon is a disk with boundary a union of geodesic segments
between $n$ vertices. 

\begin{cor}\label{cor:ngons} 
Let $S$ be an h-mapped $n$-gon. Then 
$\Area(S) \leq (n-2)\pi$. 
\end{cor}
\begin{proof} 
The exterior angle of a bent virtual vertex on a geodesic
is always less than $\pi$. 
\end{proof}

\begin{cor}\label{cor:asymzero}
Two asymptotic rays $l$ and $l'$ in $M$ or $\tilde M$ satisfy 
$d(l(t), l'(t)) \ra 0$ where $t$ is a parameter affinely related 
to the length parameter. 
\end{cor} 
\begin{proof}
Suppose not. Then there exists a sequence of pair of 
segments $(s_i, s'_i)$ meeting 
$l$ and $l'$ nearly at $\pi/2$ and the distance between 
$s_i$ and $s'_i$ goes to $\infty$. 
Therefore, there is a sequence of h-mapped quadrilaterals 
$D_i$ with sides $s_i$ and $s'_i$ and two segments $t_i$ in $l$ 
$t'_i$ and $l'$. For a sufficiently large $i$, 
there must be two points in $t_i$ and 
$t'_i$ of distance much less than the lengths of $s_i$ 
and $s'_i$ since the area of $D_i$ is bounded above by $2\pi$
and the lengths of $t_i$ and $t'_i$ go to infinity. 

Since $d(l(t), l'(t))$ is a convex function, this implies 
the corollary.
\end{proof}

\section{convex hulls in $2$-convex general hyperbolic manifolds}

Let $M$ be a $2$-convex general hyperbolic manifold 
with finitely generated fundamental group. Let 
$\core$ denote a core of $M$.


Let $\tilde M$ be the universal cover of $M$. 
Since $\core \ra M$ is a homotopy equivalence 
the subset $\tilde \core$ in $\tilde M$ of 
the inverse image of $\core$ is connected and 
is a universal cover of $\core$. 
A subset of $\tilde M$ is {\em convex} if any two 
points can be connected by a geodesic in the subset.

The {\em convex hull} $\convh(K)$ of a subset $K$ of $\tilde M$ 
is the smallest closed convex subset containing $K$. 
Since $\tilde \core$ is deck-transformation group invariant,
and the convex hull is the smallest convex subset,  
$\convh(\tilde \core)$ is deck-transformation group invariant.
Therefore, $\convh(\tilde \core)$
covers its image. We define the image as 
$\convh(\core)$, i.e., $\convh(\tilde \core)$ quotient by 
the deck-transformation group action.

Since $\core$ is a $3$-dimensional domain, 
$\convh(\core)$ is a $3$-dimensional closed set. 


\begin{prop}\label{prop:homprop} 
The convex hull $\convh(\core)$ of the compact core $\core$ of 
$M$ is homotopy equivalent to $\core$. 
\end{prop}
\begin{proof}
Let $\tilde \core$ be the inverse image of $\core$
in the universal cover $\tilde M$ of $M$. 
Then $\tilde \core$ and $\tilde M$ are both contractible 
as $M$ and $\core$ are irreducible $3$-manifolds. 

A closed curve in the convex hull $\convh(\tilde \core)$ of 
$\tilde \core$ bounds a disk since a distinguished point on
the curve can be connected by a geodesic in any other point 
of the curve. Similarly, a sphere always bounds a $3$-ball. 
Therefore, $\convh(\tilde \core)$ is contractible. 
\end{proof}





A surface is {\em pleated} if through each point of 
it passes a straight geodesic. 

Recall that the pleated-triangulated surface 
is an imbedded surface where a closed subdomain is 
a union of a locally finite collection of 
totally geodesic convex disks meeting 
each other in geodesic segments and the complementary open 
surface is pleated. 

A pleated-triangulated surface is {\em truly pleated-triangulated} 
if the triangulated part are union of totally geodesic triangles
in general position.

A {\em truly pleated-triangulated h-surface}
is a truly pleated-triangulated surface where 
each vertex of the triangles is an h-vertex.

\begin{lem}\label{lem:geohvertex}
If a geodesic in $M$ contained in $S$ 
passes through a vertex in 
the triangulated part of $S$, then the vertex is an h-vertex.
\end{lem}
\begin{proof} 
A neighborhood of a point of the triangulated part is 
stellar. If a geodesic passes through, 
the angles in both sides are greater than or equal to $\pi$:
otherwise, we can shorten the geodesic. 
Hence, the sum of the angles is greater than or 
equal to $2\pi$. 
\end{proof}

\begin{prop}\label{prop:convbdhsurf} 
Let $K$ be a deck-transformation-group invariant codimension $0$
submanifold of $\tilde M$ with $\partial K$ s-imbedded. 
Also, suppose $K$ is a subset of $\tilde M^o$.
The boundary of $\convh(K)$ can be given the structure of
a convex truly pleated-triangulated h-surface.
\end{prop} 
\begin{proof}
We will show that through each point of $\partial \convh(K)$
a geodesic in $\partial \convh(K)$ passes or 
the point is in the triangulated part and is an s-vertex 
or h-vertex. 


Let $x$ be a boundary point of $\convh(K)$:

(a) Let $x$ be a point of the manifold-interior of $\tilde M$.
Take a ball $B_\eps(x)$ in the interior for a sufficiently 
small $\eps$. Then $\convh(K) \cap B_{\eps}$ is 
the convex hull of itself. Since $B_{\eps}$ is 
isometric with a small open subset of $\hyps$, 
the ordinary convex hull theory shows that 
there exists a geodesic in the boundary of
the convex hull through $x$: If not, 
we can find a small half-open ball to decrease the 
convex hull as the side of the half-open ball cannot
meet $\partial K$ by the s-imbeddedness of $\partial K$.

Suppose that $x$ is in the topological interior of 
$\convh(K)$ but in $\partial \tilde M$. 
There exists a neighborhood of $x$ with manifold-boundary 
in $\partial \tilde M$. If $x$ is in the interior of 
an edge or a face of $\delta \tilde M$, 
then there is a geodesic through $x$ obviously. 
Suppose that $x$ is a vertex of $\delta \tilde M$. 
$x$ can be an s-vertex or a convex vertex
(see Proposition \ref{prop:s-bd2-conv}).

If $x$ is a convex vertex of $\partial \tilde M$, we can find a truncating 
totally geodesic hyperplane and a sufficiently 
small disk in it bounding 
a neighborhood of $x$ in $\tilde M$. 
Since $K$ is disjoint from the disk, 
we see that $x$ is not in the convex hull.
This is absurd. 

If $x$ is an s-vertex of $\partial \tilde M$, then $x$ is 
an s-vertex of $\partial \convh(K)$. 

Assume from now on that $x$ is a point in the topological boundary of 
$\convh(K)$. This means that $x$ is in the frontier
of the open surface $C = \partial \convh(K) - \partial \tilde M$. 

(b) Suppose $x$ is a point of the interior of a triangle $T$ in 
$\partial \tilde M$. The set $T \cap\convh(K)$ is a convex 
subset and $x$ lies in the boundary. The boundary must 
be a geodesic since we can use a small half-open ball 
to decrease the convex hull otherwise. 
Hence there is a geodesic through $x$. 

(c) Suppose that $x$ is a point of the interior of an edge
in $\partial \tilde M$. We take a small ball $B_\eps(x)$ around 
$x$, which is isometric with a ball in $\hyps$ of 
the same radius and a wedge removed. The line $l$ of the wedge 
passes through $x$. 

Let $P_1$ and $P_2$ be the totally geodesic plane extended 
in $B_\eps(x)$ from the sides of the wedge.
We denote by $P'_1$ the set $\partial B_\eps(x) \cap P_1$ 
and $P'_2$ the set $\partial B_\eps(x) \cap P_2$. 
We can form two convex subsets $L_1$ and $L_2$ in $B_\eps(x)$ 
that are the closures of the components of 
$B_\eps(x) - P_1 - P_2$ and adjacent to $P'_1$ and $P'_2$ 
respectively.  

The set $L_1 \cap \convh(K)$ is a convex subset of $L_1$
and $L_2 \cap \convh(K)$ one of $L_2$. 
The open surface $C$ may intersect $L_1$ or $L_2$ or both. 

If $C$ is disjoint from $L_1$, then 
it maybe that a one-sided neighborhood of $x$ in a triangle in 
$\partial \tilde M$ is a subset of $\convh(K)$ 
and an edge of the triangle is a geodesic through $x$.
(The side $P_1 - P_1^{\prime, o}$ of $B_\eps(x)$ 
is in $\convh(K)$.) 
Otherwise, $\convh(K)$ is contained in a convex subset 
of $B_\eps(x)$ bounded by $P_1$. In this case, 
the ordinary convexity in $\hyps$ holds and there is a geodesic 
in $\partial \convh(K)$ through $x$.


Since the same argument holds with $L_2$ as well,
we assume that $C\cap L_1$ and $C \cap L_2$ are both not empty. 

If $C \cap L_1$ is not totally geodesic for every 
sufficiently small $\eps > 0$, 
then there exists a sequence of points 
$\{x_i \in L_1 \cap \tilde M^o\}$ 
converging to $x$ and a sequence of infinitely many mutually distinct 
pleating lines 
\[\{l_i \subset \partial \convh(K) \cap \tilde M^o\}\]
so that $x_i \in l_i$. 

Therefore, by Lemma \ref{lem:pleatinglines}, 
we only have to worry about the case when $l_i$s 
end at $x$. In this case there exists a small 
neighborhood $B(x)$ of $x$ such that $C \cap L_1 \cap B$
is a cone-type set with vertex at $x$. 

By a same argument, $C\cap L_2 \cap B$ is a cone-type set also
with a vertex at $x$. In order that at $x$, the convexity 
to hold true, we see that $C\cap L_1 \cap B$ has 
to have a unique pleating geodesic and so does 
$C \cap L_2 \cap B$. They must extend each other as geodesics
passing through $x$ in order that the convex hull 
does not become less or small by truncation by 
local totally geodesic hypersurfaces.

(d) Now assume that $x$ in $\partial \convh(K)$ is 
a vertex of $\partial \tilde M$. 

Let $B_\eps(x)$ be a small neighborhood of $x$ so 
that $B_\eps(x) \cap \tilde M$ is a stellar set from $x$. 
As before $x$ is in the boundary of $C$. 

Suppose first that there are no pleating lines with a sequence of points 
on them converging to $x$. We can choose a small 
$\eps$ so that $B_\eps(x) \cap \convh(K)$ is a stellar set.

Let $M'$ be an ambient general hyperbolic manifold containing $M$
in its interior which is homeomorphic to the interior of $M$ 
as there are always such a manifold.
We claim that $x$ is an s-vertex of 
$B_\eps(x) \cap \partial \convh(K)$: 
If not, we can find a small half-open ball $B$ in $\tilde M'$ with 
a totally-geodesic side passing through $x$ with
$B^o$ disjoint from $\partial \convh(K)$. 
By stellarity, $B^o$ is disjoint from $\convh(K)$ and 
we can decrease $\convh(K)$, which is a contradiction. 
Therefore, $x$ is an s-vertex.

We assume that $B_\eps(x) \cap \convh(K)$ is not a stellar set.
Suppose now that there exists a sequence of points 
$\{x_i \in l_i\}$ converging to $x$ where $l_i$ is 
a distinct pleating line for each $i$ and does not end at $x$.
Here, $l_i$ are infinitely many.
Lemma \ref{lem:pleatinglines} shows that 
the endpoints of $l_i$ are bounded away from $x$. 
A subsequence of $l_i$ converges to a geodesic $l$ passing
through $x$. 

We see that each point $x$ of $\partial \convh(K)$ 
either has a stellar neighborhood with finite triangulations 
or has a pleated neighborhood from each of the above cases
(a), (b), (c), and (d).
We see that 
$\partial \convh(K)$ is a truly pleated-triangulated 
surface: Let $A$ be the closure of the 
set of all points in $\partial \convh(K)$ 
and in the interiors of triangles in $\partial \tilde M$.
Then $A$ is a locally finite union of totally geodesic polygons and 
segments. The complement of $A$ in $\partial K$ is pleated 
since it lies in the interior of $\tilde M$. 
Any pleated lines in $\tilde M^o$ must end at $A$ or 
is infinite. By Proposition \ref{prop:pleatinglocus}, 
the set of pleating lines ending at $A$ is isolated. 
As before, 
we see that a pleating line ending at a point of $A$ 
must be a finite length segment. 
Take a union of these finite length segments with $A$ 
to form $A'$. The infinite length pleating lines 
are in the finite union of minimal laminations 
bounded away from $A'$ by Proposition \ref{prop:pleatinglocus}.

Since the minimal laminations are bounded away from $A'$, 
we obtain a truly pleated-triangulated surface.

By Lemma \ref{lem:geohvertex}, it is an h-surface as well.
The convexity is obvious. 
\end{proof}

\begin{lem}\label{lem:pleatinglines} 
Let $l_i$, $i \in I$, be a collection of 
mutually distinct straight pleating lines 
$\partial \convh(K) - \partial \tilde M$ 
for a convex hull $\convh(K)$ of a closed subset $K$ of 
$\tilde M$ and an index set $I$. 
Suppose $x_i \in l_i$ form a sequence converging to $x$
but $x$ is not on $l_i$s. 
Then the endpoints of $l_i$s are bounded away from $x$
and a subsequence of $l_i$ converges to a line segment in the pleating 
locus containing $x$ in its interior. 
\end{lem}
\begin{proof}
Suppose that the endpoints of $l_i$ are bounded away from $x$. 
Then the second statement holds obviously.

Suppose that the endpoints $q_i$ of $l_i$
form a sequence converging to $x$. 
Then we may assume without loss of generality that 
$q_i$ lies in an arc $\alpha$ a triangle in $\partial \tilde M$.
If the arc $\alpha$ is a convex curve, we can decrease $\convh(K)$ 
further, a contradiction. Thus  $\alpha$ is a geodesic
or a point.

If $l_i$ is concurrent at a point $y$ not equal to $x$, 
then the endpoints of $l_i$ cannot be at $x$ and 
the conclusion holds. 

If all but finitely many 
$l_i$ is concurrent at $x$, then by convexity of 
$\convh(K)$, it follows that an endpoint of $l_i$ is 
$x$. This contradicts the premise. 

Suppose that no subsequence is concurrent. 
Suppose that two or more $l_i$s pass through $l$ at distinct points,
say $l_{i_j}$ for $j = 1, 2, 3$. 
Since $\convh(K)$ is convex, we can find 
a supporting line $P_j$ containing $l_{i_j}$. 
Since $P_j$ ends at a point of $l$ and an open arc in 
$l$ is a subset of $\convh(K)$, it follows that 
$l$ is a subset of $P_j$. Since $P_j$ is supporting, 
an open domain bounded by $P_j$ must contain $l_{i_l}$ and $l_{i_m}$
for $l, m \ne j$.  
However, since $P_i$s contain a common line, 
there is a $P_i$ which separates some a pair in 
$l_{i_j}, l_{i_l},$ and $l_{i_m}$.
This is a contradiction. 
\end{proof}

This proves Theorem \ref{thm:convh}.

\part{The proof of the tameness of hyperbolic $3$-manifolds}

\section{Outline of the proofs}

The strategy is as follows.
Suppose that the unique end $E$ of $M$ not associated 
with incompressible surface is not geometrically finite and 
is not tame. 
We find an exhausting sequence $M'_i$ in $M$ 
so that $M'_i$ contains neighborhoods of all tame and geometrically 
finite ends and meets the neighborhood of the infinite 
end in a compact subset.

Using the work of Freedman-Freedman \cite{FrFr}, we can
modify $M'_i$s to be compression bodies. 
Since $E$ is not geometrically finite, we can choose 
a sequence of closed geodesics $c_i \ra \infty$. 
We fix a sufficiently small Margulis constant $\eps$.
We assume without loss of generality that $c_i \subset M'_i$
and that the Margulis tubes that $M'_i$ meets are contained in 
$M'_i$. Let $\mu_i$ be the union of closed geodesics that
are in the Margulis tubes in $M_i$.
We further modify $M'_i$ so that $\partial M'_i$ is incompressible
in $M - c_1 - \cdots c_i - \mu_i$ with 
the compact core removed. The new manifolds 
$M''_i$ may form an exhausting sequence but it includes $c_i \ra \infty$. 

The manifold $N_i$ is obtained from compressing disks 
for $M'_i$. Let $A_i$ be a homotopy in $M$ between $c_i$ and 
the closed curve in $\core$. Then a compressing disk 
for $M'_i$ in $M - c_i - \core^o$ may meet $A_i$ in 
immersed circles. However, since the circles bound 
immersed disks in the compressing disk, and $c_i$ 
is not null-homotopic in $M$, we can modify $A_i$ so
that $A_i$ has one less number of components where 
$A_i$ meets the compressing disk. In this manner, 
we can find $A_i$ in $N_i$. 

Now we modify $N_i$ to $M''_i$ so that $c_i \subset M''_i$ 
and $M''_i$ are compression bodies and the boundary 
component of $M''_i$ corresponding to $E$ is incompressible 
in $M$ removed with the closed geodesics $c_1 \dots, c_n$
and contains any Margulis tubes that $M''_i$ meets.

As the boundary is incompressible, 
we take a $2$-convex hull $M_i$ of $M'_i$ using crescents
(see Part 1).
This implies that $M_i$ is a polyhedral hyperbolic space 
and hence is $\delta$-hyperbolic 
with respect to the induced path metric. 
Since $M_i$ is isotopic to $M''_i$, the homotopy 
$A_i$ between $c_i$ and a closed curve in $\core$
still exists. 

We show that we can choose a Margulis constant 
independent of $i$ and the thin part of $M_i$ are 
contained in the original Margulis tubes of $M$
and are homeomorphic to solid tori with nontrivial 
homotopy class in the original tubes. 

We take the cover $L_i$ of $M_i$ corresponding to the fundamental group of 
the fixed compact core $\core$ of $M$. Since $M_i$ is tame, the cover $L_i$ is 
shown to be tame (this is from an idea of Agol).

The core $\core$ lifts to the cover $L_i$ and can be 
considered a subset. We take a convex hull $K_i$ of 
$\core$ in $L_i$. Then $K_i$ is shown to be compact
since $L_i$ is tame and is homeomorphic to $\core$ 
if we remove parts outside a finite union of disks
disjoint from $\core$.
Since $K_i$ is homotopy equivalent to $\core$ by Proposition 
\ref{prop:homprop} in Part 2, 
the boundary component $\Sigma_i$ of $K_i$ 
corresponding to $E$ has the same genus as that of 
a boundary component of $\core$. 
$\Sigma_i$ is a ``hyperbolic surface"
(see Part 2).
Since $c_i$ is exiting an $\eps$-neighborhood of $M_i$ contains 
$c_i$, it follows that 
$p_i|\Sigma_i$ is an exiting sequence of surfaces. 

We push $\core$ inside itself so that $\core$ does not meet 
$\partial K_i$. 
We now remove the core from $K_i$ to obtain $K_i - \core^o$.
We can find a simple closed curve $\alpha$ in 
$\Sigma_i$ compressible in $K_i$. 
We realize $\alpha$ by a closed geodesic 
$\alpha^*$ in $K_i - \core^o$. Since $\alpha$ cannot 
be realized by a straight geodesic or even a quasi-geodesic, 
$\alpha^*$ must 
meet $\partial \core$. Fixing a base point $x^*$ on 
$\alpha$, we choose relative geodesics on $\Sigma_i$
which triangulates $\Sigma_i$ with $\alpha$.
Homotopy $x^*$ to a base point $y^*$ on $\alpha^*$ 
and geodesics to have endpoints in $y^*$.
We obtain a simplicial 
hyperbolic surface $T_i$ meeting $\partial \core$.

Then by compactness of bounded simplicial
hyperbolic surfaces, infinitely many immersed $T_i$s are isotopic
in $M - \core^o$ and hence infinitely many of $p_i|\Sigma_i$ are 
isotopic. Since $p_i|\Sigma_i$ are exiting
and are isotopic in $M - \core^o$,
this implies that the end $E$ is topologically tame.

We remark that the most important new idea of Agol, originally 
from Freedman, is the use of tameness
of the covering spaces of tame manifolds and taking convex hulls 
in the covers. These were the approaches that we adopted in
this paper. 



\section{The Proof of Theorem A} 


We will now choose the compact core more carefully so that 
$\partial \core$ is s-imbedded. 

We can choose incompressible closed 
surfaces $F_i$ associated with incompressible 
ends $E_i$ to be strictly s-imbedded by Theorem C and disjoint 
from one another (see Remark \ref{rem:svertices}).
We choose a number of closed geodesics in 
$E_i$ and choose a mutually disjoint submanifold homeomorphic to 
$F_i \times I$ disjoint and between these curves for each $i$.
Then by Theorem C, we find a mutually disjoint collection of
manifolds in the respective neighborhoods of $E_i$
between these curves whose boundary components 
are strictly s-imbedded.


Essentially $\partial \core$ can be considered 
a regular neighborhood of the union of the essential surfaces 
$F_1, \dots, F_n$ and a number of arcs connecting them 
in some manner. 

We choose each of the arcs to be the shortest 
path in $M$ among the arcs connecting 
the surfaces $F_i$s with the respectively given homotopy classes.
Their endpoints must be in the interior of 
an edge of a triangle. 
By perturbing $F_i$ if necessary, we may 
assume that they are all disjoint geodesics. 
We first take thin regular neighborhoods of $F_i$s 
which are triangulated. 
We take thin regular neighborhoods of the geodesics
which are triangulated and all of whose vertices lie in
$F_i$s. 

We take the union of the regular neighborhood of these 
geodesics with those of $F_i$s to be our core $\core$. 
we may assume that $\partial \core$ is strictly s-imbedded as well.
(We may need to modify a bit where the neighborhoods meet.)
As stated above, we choose $\core$ to be in the interior 
$M^o$. Obviously, if necessary, we can push $\core$ inward 
itself without violating strict s-imbeddedness of $\partial \core$. 


Let $M$ be as in the introduction, and let $U_1, \dots, U_n$ 
be mutually disjoint neighborhoods of incompressible ends
$E_1, \dots, E_n$. Suppose that the end $E$ is a geometrically 
infinite but not geometrically tame. 

Let $\hat M$ be the $2$-convex hull of 
$M$ with $U_1, \dots, U_n$ removed.
The boundary components $F_1, \dots, F_n$ 
corresponding to $U_1, \dots, U_n$ 
of $\hat M$ are s-imbedded respectively.

Let $M'_i$ be an exhaustion of $M$
by compact submanifolds containing $\hat M$. 
We extend $M'_i$ by taking a union with $U_1, \dots, U_n$
so that $M'_i$ meets neighborhoods of $E$ in 
compact subsets or in the empty set. 
We assume that $M'_i$ contains 
the boundary components $F_1, \dots, F_n$ and 
contains the core $\core$ of $M$ always.

We fix a small Margulis constant $\eps_M >0$.

\begin{prop}[Freedman-Freedman, Ohshika]
\label{prop:compressiong}
We can obtain a new exhaustion $M'_i$ so that 
each $M'_i$ is homeomorphic 
to a compression-body with incompressible
boundary components removed. $M'_i$ contains 
the Margulis tubes that $M'_i$ meets
by taking union with these.
\end{prop}
\begin{proof} 
We assume that each $M'_i$ includes any Margulis tubes it meets. 
We essentially follow Theorem 2 of Freedman-Freedman \cite{FrFr}. 
$\partial M$ has no incompressible closed surface other than 
ones parallel to $F_1, \dots, F_n$. 
Hence, we can compress the boundary component 
$\partial M'_i$ until we obtain a union of 
balls and manifolds homeomorphic to 
$F_i$ times an interval. 
An exterior disk compressions add a disk times $I$ to the compressed 
manifolds from $M'_i$ but interior disk compressions remove 
a disk times $I$ from the manifolds.
The exterior disk may meet Margulis tubes outside 
$M'_i$. If the disk meets a Margulis tube 
at an essential disk, then we include the Margulis tube.
If not, we push the disk off the Margulis tube. 
This operation amounts to adding $1$-handles to 
the manifolds.

We recover our loss to $M'_i$ by interior disk compressions
by attaching $1$-handles each time we do 
interior compressions. The core arcs of $1$-handles 
may meet the exterior compression disk many times. 
We add a small neighborhood of the cores first. 
Then we isotopy to make it larger and larger to recover 
the loss due to interior disk compression
while fixing the Margulis tubes outside $M_i$. 
(This may push the exterior compression disks.)
We also recover all the Margulis tubes originally in $M'_i$.

From the surface times the interval components, we 
add all the $1$-handles to obtain the desired compression body.

\end{proof} 

\begin{rem}\label{rem:intext}
We do interior compressions and 
then exterior compressions. This is sufficient
to obtain the union of $3$-balls and $F_i$ times 
the intervals. The reason is that we can isotopy 
any interior compressing curve away from the traces of 
exterior compression disks.
\end{rem}


Since $E$ is geometrically infinite, there exists 
a sequence of closed geodesics $c_i$ tending to $E$
by Bonahon \cite{Bonahon}.
We assume that $c_i \subset M'_i$ for each $i$
since $M'_i$ is exhausting. Let $\calC_i$ denote the union of 
$c_1, \dots, c_i$.
We assume without loss of generality that 
$M'_i$ has a free-homotopy $A_i$ between $c_i$ and a 
closed curve in $\core$.   
Let $\mu_i$ be the union of the simple closed geodesics 
in the Margulis tubes that $M'_i$ contains. 

\begin{rem}\label{rem:diskslide}
We know that given a surface there exists 
a finite maximal collection of exterior essential compressing 
disk so that any other exterior essential compressing 
disks can be pushed inside the regular neighborhood of their union.
This is essentially the uniqueness of the compression body.  
(See Theorem 1 Chapter 5 of McCullough \cite{McNote}.)
\end{rem}

If we remove the interior $\core^o$ of the core from $M'_i$, 
and compress the boundary 
$\partial M'_i$ in $M - \core^o - \calC_i - \mu_i$ 
and remove resulting cells to obtain a manifold 
with incompressible boundary containing $\partial \core$. 
Then we join the result with $\core$, $\calC_i$ and $\mu_i$.
Let us call the resulting $3$-manifold $N_i$. 
Note that $N_i$ need not be a compression body 
and may not form an exhausting sequence. 
However, $N_i$ contains $\core, \calC_i, \mu_i$
for each $i$.

The exterior compressing disk of $M'_i$ 
may meet the Margulis tubes outside $M'_i$. 
We may assume that $N_i$ meets these Margulis 
tubes in meridian disks times intervals
by following Lemma \ref{lem:mtubes}.

\begin{lem}\label{lem:mtubes}
A disk with boundary outside the union of 
Margulis tubes may be isotopied with the boundary 
of the disk fixed so that 
the intersection is the union of meridian disks. 
\end{lem}
\begin{proof}
First, we perturb to obtain transversal intersections. 
A disk may meet the Margulis tubes in 
a union of planar surfaces. 
The boundary of the Margulis tube meets the disk 
in a union $C$ of circles. If an innermost component 
is outside the tube, then since the boundary tube is 
incompressible in $M$ with the interior of 
the Margulis tubes removed, it follows that we 
can isotopy it inside. This means that 
the innermost components are disks. 

If an innermost component of $C$ is a circle bounding 
a disk in the boundary of the Margulis tube, 
then we can isotopy the bounded disk away from the tube. 
Now, each component of $C$ is a meridian circle. 
\end{proof}

The manifold $N_i$ is obtained from compressing disks 
for $M'_i$. Let $A_i$ be a homotopy in $M'_i$ between $c_i$ and 
a closed curve in $\core$. Then a compressing disk 
for the sequence of manifolds obtained from $M'_i$ 
by disk-compressions
in $M - \core^o - \calC_i - \mu_i$ may meet $A_i$ in 
immersed circles. However, since the circles bound 
immersed disks in the compressing disk, and $c_i$ 
is not null-homotopic in $M$, we can modify $A_i$ so
that $A_i$ has one less number of components where 
$A_i$ meets the compressing disk. In this manner, 
we can find $A_i$ in $N_i$.

We do the following steps:
\begin{itemize} 
\item[(1)] We find a maximal collection of 
essential interior compressing disks for $N_i$
and do disk compressions. By incompressibility,
the disk must meet one of $\core, \calC_i, \mu_i$ essentially.
We call $N^1_i$ the component of the result containing $\partial \core$. 
\item[(2)] We find a maximal collection of 
essential exterior compressing disk for the result of (1) 
and do disk compressions. We call $N^2_i$ the component 
containing $\partial \core$. 
\item[(3)] We add $1$-handles lost in the step (1).
We call the result $N^3_i$. 
\end{itemize} 


Clearly $N^3_i$ includes $N_i$.

\begin{prop}\label{prop:incompressible2} 
The submanifold $N^3_i$ is homeomorphic to a compression body 
and is contained in a compression body $M''_i$ whose 
boundary component is incompressible in 
$M - \core^o - \calC_i - \mu_i$. A Margulis tube is 
either contained in $M''_i$ or the tube meets $M''_i$ 
in meridian disk times intervals. 
\end{prop}
\begin{proof} 
The fact that $N^3_i$ is a compression body follows as before. 
Using the fact that 
$N^3_i$ is contained in some compression body $M'_j$ for a large $j$,
we will now show that $N^3_i$ is contained in
a compression body $M''_i$ with the above property.

By construction of $N^3_i$, it follows that 
an interior compression disk is equivalent to one 
of the disks of the $1$-handles by disk sliding.
Therefore, each interior compressing disk must 
intersect at least one of $\core, \calC_i, \mu_i$ essentially.
Hence $\partial N^3_i$ has no interior compression 
disk in $M - \core^o - \calC_i - \mu_i$. 

There could be an exterior compression disk for $N^3_i$. 
We take a maximal mutually disjoint 
family $D_1, \dots, D_n$ of them where no two of 
$\partial D_i$ are parallel.
We choose $j$ sufficiently 
large so that a compression body $M'_j$ includes all of them
as $M'_i$s form an exhausting sequence.

We find a $3$-manifold $X$ isotopic to $N^3_i$ in $M'_j$:
$M'_j$ decomposes into a union of cells or submanifolds homeomorphic to
$F_l$ times intervals by a maximal family of interior compression 
disks. We suppose that no two of the disks are parallel. 
We call $D$ the union of these disks.  

We choose $X$ so that $X$ is a thin regular neighborhood of 
a $1$-complex with the unique vertex in a fixed base cell of $M'_j$.
We fix a handle decomposition of $X$ corresponding to the $1$-complex 
structure. We define the complexity of the imbedding of $X$ in $M'_j$ 
by the number of components of $X \cap D$.
We choose $X$ with minimal complexity. 
We will now put all things in general positions.
For each disk $D_k$, we first get rid of any closed circles 
by the innermost circle arguments. 
We will find an edgemost arc if 
$\partial D_k$ meets $D$ bounding a component of 
$D_k -D$. Then there must be a handle of 
$X$ following the arc in $\partial D_k$. 
This handle can be isotopied away,
and then using the innermost circle argument again 
if necessary, we can reduce the complexity. 
Therefore, it follows that 
each $D_k$ is in the fixed base cell of $M'_j$. 
Also, a handle where $D_k$ passes essentially
must lie in the base cell also.

We look at the handles in the base cell and 
the disks. The union of the handles and the ball around 
the basepoint is a handle body.
Our disks $D_1, \dots, D_n$ are in the cell. 

From Corollary 3.5 or 3.6 of Scharlemann-Thomson \cite{ScTh}, 
we see that there exists an unknotted cycle in the $1$-complex 
or the $1$-complex has a separating sphere. 
In the first case, we cancel the corresponding cycle 
by an exterior compression. In the later case, 
a sphere bounds a ball which we add to $X$, i.e., we engulf it. 
We do the corresponding topological operations to 
$N^3_i$. Note that $X$ and $N^3_i$ are still compression bodies. 
In both cases, we can reduce the genus of 
the boundary of the handle body $X$ or $N^3_i$.
We do this operations until there are no more 
exterior compressing disks. (The genus complexity shows 
that the process terminates.)
We let the final result to be $M''_i$.

The imbedding $\partial M''_i \ra M - M^{\prime \prime, o}_i$ 
is incompressible since the boundary of any exterior compressing disk 
for $\partial M''_i$ can be made to avoid the traces of 
handle-canceling exterior compressions which are pairs of disks
or the disks from the $3$-ball engulfing.
Therefore these must be exterior compressing disks for 
$N^3_i$. 

The imbedding $\partial M''_i \ra M^{\prime \prime, o}_i 
- \core^o - \calC_i - \mu_i$ 
is incompressible since the boundary of any interior compressing disk
can be made to avoid the traces of handle-canceling 
exterior compressions.

The statement about Margulis tubes follows by Lemma \ref{lem:mtubes}.
\end{proof}


We will now modify $M''_i$ by crescent-isotopy.

\begin{lem}\label{lem:crescentcore} 
A secondary highest-level crescent of 
$\tilde \Sigma$ does not meet the interior of $\tilde \core$. 
\end{lem}
\begin{proof} 
If not, then $\tilde \core^o$ meets $I_\calR$ for a secondary 
highest-level crescent $\calR$. We may assume that $\calR$ is 
compact by an approximation inside. 
Again find a Morse function by totally geodesic planes 
parallel to $I_\calR$. $\core \cap \calR$ has a maximum 
inside as $\core \cap \calR$ is compact. 
But at the maximum point, a totally geodesic plane 
bounds a local half open ball disjoint from $\core^o$. 
This contradicts the s-imbeddedness of $\partial \tilde \core$. 
\end{proof}

We define $M_i$ to be the $2$-convex hull of $M''_i$.
The boundary components of $M_i$ are s-imbedded.
Let $\core'$ be the core obtained from $\core$ by pushing $\partial \core$
inside $\core$ by an $\eps$-amount.
Note that during the crescent-isotopy
$\core'$ is not touched by the interior of secondary highest 
level crescents since 
$\partial \core$ is s-imbedded by Lemma \ref{lem:crescentcore}.

Define $M_i^\eps$ be the regular $\eps$-neighborhood of $M_i$.
$\calC_i, \mu_i \subset M_i^\eps$ 
by Proposition \ref{prop:geodesicavoid} in Part 1.
We may assume $A_i$ is in $M_i^\eps$ since we isotopied $M''_i$ 
to obtained $M_i$.

The universal cover $\tilde M_i$ of $M_i$ 
with the universal covering map $p_i$
is an $M_{-1}$-space and is $\delta$-hyperbolic. 

We define the {\em thin part} of $M_i$ as the subset of $M_i$ 
where the injectivity radius is $\leq \eps$. 
Since $\tilde M_i$ is a uniquely geodesic,  
through each point of the thin part of $M_i$
there exists a closed curve of length $\leq \eps$
which is not null-homotopic in $M_i$.

Let $\gamma$ be a closed curve of length $\leq \eps$ 
which is not null-homotopic in $M_i$. 
Then if $\gamma$ has nontrivial holonomy, 
then $\gamma$ lies in a Margulis tube of $M$ 
which either is in $M_i$ or disjoint from $M_i$. 

Suppose that 
$\gamma$ is null-homotopic in $M$. 
Then $\gamma$ bounds an immersed disk $D$ in $M$. 
The diameter of $D$ is $\leq \eps$.

Suppose that $D$ cannot be isotopied into $M_i$.
By incompressibility of $\partial M_i$, 
$D$ must meet $c_i$ or $\mu_i$ or $\core$. 
Also, $D/M_i$ must be nonempty.
There must be an innermost disk $D'$ such that  
$\partial D'$ maps to $\partial M_i$ and 
$D'$ maps into $M_i$ and 
meets $\mu_i$, $\core$ or $\calC_i$
and a component of $D/\partial M_i$ 
adjacent to $D'$ lies outside $M_i$.

If $D'$ meets $\core$, then the diameter of $\partial D'$ 
is not so small, and hence that of $D$. 
If $D'$ meets $\calC_i$, then since $A_i$ is in $M_i$ 
$D'$ cannot be bounded outside by a component 
in $M /M_i^o$. 

Suppose that an innermost disk $D'$ meet $\mu_i$. 
Since $\partial D'$ is very close to $\mu_i$
due to its size, and 
the distance from $\mu_i$ to $\partial M''_i$ 
is bounded below by a certain constant, 
it follows that the boundary of $D'$ lies 
in the union of $I$-parts of some crescents
or its perturbed images
obtained during the crescent-isotopies. 
Since the length of components of $\mu_i$ is short, 
we see that the the $I$-parts meeting $\partial D'$
would extend for long lengths along the geodesics
near a component of $\tilde \mu_i$. 
Therefore, we see that at the last stage of
the isotopies, we have the inverse image of 
torus bounding a component of $\mu_i$. 
Since our crescent moves are isotopies and $\partial M'_i$ is 
not homeomorphic to a torus, this is a contradiction. 

We conclude that $\gamma$ bounds a disk in $M_i$. 
Since $\gamma$ is not null-homotopic in $M_i$, 
this is a contradiction. 
Therefore, the thin part of $M_i$ is in the intersection of 
the Margulis tubes of $M$ with $M_i$. 

Note here that the Margulis constant $\eps >0$ could be 
chosen independent of $i$ as the above argument 
passes through once $\eps$ is sufficiently small 
regardless of $i$.

By above discussions, it follows that any $\eps$-short closed 
curve in $M_i$ is a multiple of the simple closed geodesic 
in a  Margulis tube. We may assume that given a component of 
the thin part of $M_i$,
an $\eps$-short closed curve of a fixed homotopy class passes  
through each point of the component.
Therefore, components of thin parts are solid tori in Margulis 
tubes. 

During the crescent move, the shortest closed 
geodesic in the Margulis tube may go outside 
particularly during the convex perturbations.
But there are short closed curves in the result that homotopic 
to the closed geodesic.
Therefore the thin part are union of solid tori 
parallel to the shortest geodesics.

\begin{prop}\label{prop:Margulis}
The $\eps$-thin part of $M_i$ are homeomorphic to a disjoint union of 
solid tori in Margulis tubes in $M$ parallel to 
the shortest geodesics in the respective tubes.
Furthermore, we can choose $\eps >0$ independent of $i$. 
\qed
\end{prop}

Assume without loss of generality that 
we have an inclusion map $i:\core \ra M_i$ for each $i$. 
Let $L_i$ be $\tilde M_i/i_*\pi_1(\core)$
with the covering map $p_i:L_i \ra M_i$. 
$L_i$ has ends corresponding to $F_1, \dots F_n$, 
and another end $E$ corresponding to $E$. (We abused notation 
a little here.)

\begin{prop}\label{prop:compK_i}
The convex hull $K_i$ of $\core$ in $L_i$ is compact.
\end{prop} 
\begin{proof}
Since $M_i$ is tame, its cover $L_i$ is tame. 
For any compact set, we can find a compact core 
of $L_i$ containing it. By choosing a large 
compact subset of $M_i$, we obtain a compact 
core $\core'$ of $L_i$ containing it which is obtained 
as the closure of the appropriate component of $L_i$ with 
a finite number of disks removed.

Certainly $\core$ is a subset of it. 
The disks lifts to disks in the universal cover $\tilde L_i$ 
of $L_i$. They bound the universal cover $\tilde \core'$ of $\core'$. 

The convex hull of a disk is a compact subset of 
$\tilde L_i$ since the convex hull of a compact 
subset is compact in the universal cover.
Therefore, the convex hull of $\tilde \core'$ is in the union 
of $\tilde \core'$ and deck-transformation images of 
finitely many compact sets 
that are convex hulls of the boundary disks. 
Therefore, the convex hull of $\core$ itself is compact
being a subset of the union of a compact set $\core'$ and 
finitely many compact sets.  
\end{proof}

Since $\core$ is homotopy equivalent to $L_i$, 
there exists a convex hull $K_i$ of $\core$ in $L_i$ homotopy 
equivalent to $\core$ by Proposition \ref{prop:homprop} in Part 2.
Obviously, $K_i$ contains $F_1, \dots, F_n$. 
Let $\Sigma_i$ be the unique boundary component of 
$\partial K_i$ associated with $E$. 

Any closed geodesic homotopic to 
a closed curve in $\core$ in $L_i$ is contained in $K_i$: 
If not, we can find an h-imbedded annulus $B_i$ 
with boundary consisting of  
the closed geodesic and a closed curve on $\partial K$, 
intrinsically geodesic,  
where the interior angles in $B_i$ are always greater than or equal to 
$\pi$. Such an annulus clearly cannot exist.  

We can find a quasi-geodesic $c'_i$ $\eps$-close to $c_i$ 
in $M_i$. Since $c'_i$ is homotopic to a closed curve in 
$\core$ by a homotopy $A'_i$ in $M_i$ modified from $A_i$,
we have that $c'_i \subset L_i$.
Choose a geodesic representative $c''_i$ in $L_i$
which is again arbitrarily close to $c'_i$ and hence to $c_i$.
Therefore $c''_i \subset K_i$ for each $i$.

Since $\Sigma_i$ is a truly pleated-triangulated convex h-surface,
the intrinsic metric in the pleated part is a Riemannian 
hyperbolic ones. Thus, $\Sigma_i$ carries a triangulated 
h-surface structure intrinsically. 
Since $\Sigma_i$ is intrinsically an h-imbedded surface, 
$p|\Sigma_i$ is one also and hence form an exiting sequence in $E$. 

More precisely, the parts of boundary of the image of $K_i$ 
form an exiting sequence in $E$. Any part of the boundary of 
the image of $K_i$ is in the image of $\Sigma_i$. 
Hence, there exists an exiting sequence of parts 
of $\Sigma_i$. By the uniform boundedness of $\Sigma_i$, 
it follows that $\Sigma_i$ form an exiting sequence in $E$. 

\begin{rem}\label{rem:e-thinh-surf}
The uniform nature of the Margulis constant plays 
a role here. Any $\eps$-thin part of an h-immersed surface 
must be inside a Margulis tube in $L_i$ and by 
incompressibility the thin parts are homeomorphic 
to essential annuli. Since $\Sigma_i$ is incompressible 
in $L_i - \core^o$, we see that the thin part of $\Sigma_i$ 
is a union of essential annuli which are not homotopic 
to each other. Thus, outside the Margulis tubes, 
$\Sigma_i$s have bounded diameter independent of $i$. 
\end{rem}

\section{The Proof of Theorem B}

We recall that $\core$ was pushed inside itself somewhat so that 
$\Sigma_i$ and $\partial K_i$ does not meet $\core$. 
In $K_i$, we may remove $\core^o$ and 
we obtain a compact submanifold $Q_i$ of 
codimension $0$ bounded by s-surfaces
including $\Sigma_i$ since $\partial \core$ is s-imbedded. 
$Q_i$ is $2$-convex and is Gromov-hyperbolic. 
$\Sigma_i$ is incompressible in $Q_i$
since any compressing disk of $\Sigma_i$ not meeting the core 
would reduce the genus of $\Sigma_i$ but the genus of 
$\Sigma_i$ is the same as that of the component of $\partial \core$
corresponding to the end $E$ since $K_i$ is homotopy equivalent 
to $\core$.

As $K_i$ is homeomorphic to a compression body,
we choose a compressing curve $\alpha$ in $\Sigma_i$.
Then $\alpha$ bounds a disk $D$ in $K_i$
and the core $\core$ must meet $D$ in its interior. 
Let $\hat \alpha$ be the geodesic realization of 
$\alpha$ in $K_i - \core^o$, which must be in $K_i$. 

If $\hat \alpha$ does not meet $\partial \core$, 
then it maps to a geodesic in $M$. This is absurd since 
the holonomy of $\alpha$ is the identity. 
Therefore, $\hat \alpha$ meets $\partial \core$. 

We form a triangulation of $\Sigma_i$ with the only vertex $p$
at a point of $\alpha$ and including $\alpha$ as an edge. 
Then choosing a vertex $\hat p$ in $\hat \alpha$ 
and a path from $p$ to $\hat p$, we isotopy each edge of 
the triangulation to a geodesic loop in $K_i - \core^o$ 
based at $\hat p$. Each triangle is isotopied to $2$-A-simplex 
spanned by new geodesic edges. The resulting 
surface $T_i$ is an h-surface since each of the triangles 
is a $2$-A-simplex and a geodesic passes through each 
point of the $1$-complex.

Each surface $q_i: T_i \ra M - \core^o$
has the same genus and is homotopic to 
$p_i| \Sigma_i$ in $M - \core^o$. 
Since they are h-imbedded, and meet $\partial \core$, 
they are in a bounded neighborhood of $\core$ by the boundedness of h-imbedded 
surfaces. They form a pre-compact sequence. 
Thus infinitely many of $q_i|T_i$ are isotopic 
in $M - \core^o$. Therefore, infinitely many of $p_i|\Sigma_i$ 
are isotopic in $M-\core^o$. 
Since $\Sigma_i$ bounds larger and larger domains in a cover of $M$ 
and $\Sigma_i$ projects to a surface far from $\core$, 
the above fact shows that our end $E$ is 
tame as shown by Thurston \cite{Thnote}.
(This is essentially the argument of Souto \cite{Souto}.)


\bibliographystyle{plain}

\end{document}